\documentclass[11pt,leqno]{article}
\usepackage{setspace}
\usepackage{graphicx}
\usepackage{amssymb}
\usepackage{epstopdf}
\title{Shannon entropy for stationary processes and dynamical systems}

\author{ By D. Hamdan,  $\hskip 0,2 cm $W. Parry\dag $\hskip 0,2 cm $  and  $\hskip 0,2 cm $J.-P. Thouvenot}
\begin{document}
\maketitle

\begin{abstract}{We consider stationary ergodic processes indexed by $\mathbb Z$ or $\mathbb Z^n$ whose  finite dimensional marginals have laws which are absolutely continuous with respect to Lebesgue measure.\\
We define an entropy theory for these continuous processes, prove an analog of the Shannon Breiman Mac Millan theorem and study more precisely the particular example of Gaussian processes.}\end{abstract}

 \section{Introduction}
    In [17] Shannon introduced a general theory of entropy designed to quantify the rate at which information is produced through the evolution of a stationary stochastic process. He gave two definitions: one for random variables assuming discrete values and the other for random variables which assume  real values or, more generally, values in $\mathbb R^n$. However these definitions represent but two aspects of a single notion.\\
     In 1958 Kolmogorov adapted Shannon's discrete version in his definition of entropy for dynamical systems which enabled him to solve an important open problem concerning the classification of measure-preserving transformations with continuous (in fact Lebesgue) spectrum. Shortly afterwards Sinai modified and improved Kolmogorov's definition so that today one speaks of the Kolmogorov-Sinai, K.S., invariant.\\
     The K.S. entropy is also an invariant for stationary stochastic processes in as much as they may be represented  as measure preserving transformations. However these processes frequently have infinite entropy particularly when the random variables take their values in a non-discrete space.\\
     For example all stationary Gaussian processes with absolutely continuous spectrum have the same infinite K.S. entropy and indeed they are all measure theoretically isomorphic to each other (Ornstein [9]).\\
     For this reason we feel that a modification of the traditional definition of isomorphism should be considered which distinguishes between various processes with infinite K.S. entropy. This paper should be regarded as a move toward this end in that we produce an $^{\ll}$invariant$^{\gg}$ based on Shannon's second version of entropy (designed for continuous valued random variables).\\
     Shannon's entropy was extensively investigated, especially by the Russian school, in the late 1950's  and a thorough account appears in Pinsker's book [14]. One of our main purposes is to clarify and extend known results in this area.\\
     The invariant investigated here ( naturally referred to as Shannon entropy ) is simply a normalised limit of entropy for $\mathbb R^n$ valued random variables, which we compute for all stationary Gaussian processes in terms of their spectral measures. In general this entropy is finite even when the K.S. entropy is infinite. We  show how this entropy changes when a stationary Gaussian process is subjected to a linear transformation.\\
     This is a topic closely related to the work of Wiener [20] and Kolmogorov [7] on linear prediction theory. It is well known that stationary Gaussian processes may be regarded as non-linear extensions of stationary linear sequences in Hilbert space (i.e. $^\ll$wide sense$^\gg$  processes in the language of Doob [2] ) and entropy theory may be regarded as the prediction theory of non-linear processes. Wiener was particularly interested in non-linear prediction and in [21] he attempted to prove that all stationary processes satisfying certain mild conditions are isomorphic to independent Gaussian processes. As shown by Rosenblatt [15] his proof was flawed and it is interesting to note that his mistake, which concerned the behaviour of decreasing sequences of sigma algebras, was repeated by Kolmogorov (as shown By Rokhlin) in a different context some years later. We now know through the work of Ornstein and his co-workers that Wiener's claim is actually false.\\
     We shall be concerned almost exclusively with stationary stochastic processes $X=(X_n)$ for which the distribution measure of $(X_0,...,X_n)$ is absolutely continuous with respect to $n+1$ dimensional Lebesgue measure for all $n=0,1,...$ Our first theorem shows that for any ergodic measure-preserving transformation $T$ of a probability space there exist functions $F$ such that $(F\circ T^n)$ satifies the above condition.
We prove thereafter that an entropy $^\ll$ \`a la Kolmogorov$^\gg$ can  be defined using the continuous entropy definition of Shannon. We prove that  these averages of entropies always converge. The theory extends naturally to $\mathbb Z^n$ actions. In this framework, we prove a Shannon Mac Millan Breiman type pointwise theorem. However,  if the Kolmogorov-Sinai entropy of the transformation is finite, then the limit above is always $-\infty$.  In general, this limit is majorized by ${1\over 2} log(2\pi)$ plus  one half of the variance of the observable. In particular, the preceding inequality reduces to equality if and only if the process is Gaussian independent ( This generalizes to processes a result of Shannon for  the case of random variables).  We give also a similar characterization for a stationary process to be Markovian, and,  more generally, to be with memory $p$. The theory applies naturally to Gaussian processes, for which,   we give a closed formula for this continuous entropy. It turns out that, in the Gaussian Markovian case, this entropy determines the process, and in the Gaussian case, when  finite it determines the process up to unilateral isomorphism. \\
   We give also some  relationships  with the rate of entropy of Pinsker and with the rate of generation of information as well.\\
  William Parry, our co-author,   died  August 20-th 2006,  at the time where these notes were  completed.
     \section{Absolutely continuous processes based on an ergodic system} 
We begin with\\ 
{\bf{Definition 2.1:}}\\
\textit{A real valued discrete time process $X=(X_n)_{n\in I}$, indexed by a countable set $I$,  is said to be absolutely continuous if for every finite subset  $K$ of $I$, with cardinality $\mid K\mid$, the joint distribution of $(X_n)_{n\in K}$  is absolutely continuous with respect to the $\mid K\mid $-dimensional Lebesgue measure.}\\
We shall only be using $I=\mathbb Z$, or $I=\mathbb Z^2$.\\
In this section we prove the following:\\
{\bf{Proposition 2.2:}}\\
\textit{Let $(\Omega,T,\mu)$ be an invertible ergodic dynamical system.  Then there exists $F\in L^2(\mu)$ such that the process $(F\circ T^n)$ is absolutely continuous.}\\
{\bf{Proof:}}  Let $(Y,S,\gamma)$ be the independent gaussian dynamical system: $Y:=\mathbb R^\mathbb Z$, $S$ the shift transformation: $(Sy)_j=y_{j+1}$, for $y\in Y,j\in \mathbb Z$, and $\gamma$ the product measure of the measures $\gamma_j,j\in \mathbb Z$, where, for every $j$,  $\gamma_j=(2\pi)^{-{1\over 2}}exp(-{1\over 2}x^2)dl(x)$,  $l$ being the Lebesgue measure on $\mathbb R$. Let $Y_n:Y\rightarrow \mathbb R$ be the projection onto the $n$'th coordinate. \\ According to Dye's Theorem, [3],  there exists an integer-valued measurable function $\tau:Y \rightarrow \mathbb Z$ such that if $S_1:Y\rightarrow Y$ is the transformation defined by   $S_1(y)=S^{\tau(y)}(y)$, for $\gamma$-almost all $y\in Y$, then the two dynamical systems $(\Omega,T,\mu)$ and $(Y,S_1,\gamma)$ are isomorphic. Let $\theta:\Omega \rightarrow Y$ be a map giving the isomorphism: $\theta \circ T=S_1\circ \theta$ and $\gamma=\mu\circ \theta^{-1}$. Set $F=Y_0\circ \theta$ and $Z_n=Y_0\circ S_1^n$, for $n\in \mathbb Z$, so that $(F,F\circ T,...,F\circ T^{n-1})$ and $(Z_0,Z_1,...Z_{n-1})$ have the same law, say $\alpha_n$. We show now that $\alpha_n$ is absolutely continuous with respect to Lebesgue measure and this  ends the proof. To do this,  for every  $\bar{k}=(k_1,...,k_{n-1})\in \mathbb Z^{n-1}$,  put 
  $$E_{\bar{k}}=\{\tau=k_1,\tau\circ S^{k_1}=k_2,\tau\circ S^{k_1+k_2}=k_3,...,\tau\circ S^{k_1+...+k_{n-2}}=k_{n-1}\},$$
  and let
${\cal F}_{\bar{k}}$ denote the sigma-algebra generated by $Y_0,Y_{k_1},...,Y_{k_1+...+k_{n-1}}$. In view of the definition of $(Y,S,\gamma)$,   straightforeward computations show then that $\alpha_n$ is absolutely continuous with respect to Lebesgue measure and has the density $ g:=\sum_{\bar{k}\in Z^n} g_{\bar{k}}$, where 
$$g_{\bar{k}}(y_0,...,y_{n-1})=E_\gamma [1_{E_{\bar{k}}}\mid {\cal F}_{\bar{k}}](y_0,y_1,...,y_{n-1})\times (2\pi)^{-{n\over 2}}exp(-{1\over 2}(y_0^2+...+y_{n-1}^2)),$$
and $E_\gamma [1_{E_{\bar{k}}}\mid {\cal F}_{\bar{k}}]$  denotes the conditional expectation with respect to $\gamma$ of $1_{E_{\bar{k}}}$ given the sigma algebra $ {\cal F}_{\bar{k}}.\square$\\
 {\bf{Remark 2.3:}} \\
 (1)  The function $F$ in Lemma 2 can be taken ( as the proof shows ) in the intersection of 
$\{L^p(\mu):p\ge 1\}$.\\
(2)  We can show that there is $F$  such that, for every $n$, the law of $(F,F\circ T,...,F\circ T^{n-1})$ is equivalent to the $n$-dimensional Lebesgue measure.\\
In the same way we have\\
{\bf{Proposition 2.4:}}\\
\textit{If $T$ and $S$ are measure preserving transformations which commute on a probability space $(\Omega,{\cal F},\mu)$, such that the joint action is ergodic then there exists $F\in L^2(\mu)$ such that the process $(F\circ T^m\circ S^n)_{(m,n)\in \mathbb Z^2}$ is absolutely continuous.}

\section{Notations, a few prerequisites and a lemma}
\subsection{Conditional entropy of probability measures}
In this subsection we recall various definitions attached to Shannon entropy [14
].\\
\textbf{Definition 3.1 :}\\
\textit{Let $P$ and $Q$ be  probability measures defined on the same measurable space $(\Omega,{\cal F})$. Let $\cal P$ be the set of all finite measurable partitions of $\Omega$. If $\Pi$ is in  $\cal P $   let
\begin{eqnarray}S_\Pi (P\mid Q):=\sum_{E\in \Pi} P(E)log({P(E)\over Q(E)}).\end{eqnarray}
Then for $\Pi_1\in \cal P $, $\Pi_1$  finer than $\Pi$ implies 
\begin{equation}
S_\Pi (P\mid Q)\le S_{\Pi_1} (P\mid Q).
\end{equation}
  The entropy $H_{Q}(P)$ of $P$ with respect to $Q$ is defined by
  \begin{eqnarray}
  H_Q(P):=\sup_{\Pi \in \cal P} S_\Pi(P\mid Q).\end{eqnarray}  
}

We list without proofs some results which we are going to use.\\

{\bf{Theorem A (Gelfand, Yaglom, Perez, [14]):}}
\textit{Let  $P,Q$ be two probabiliy measures on the measurable space $(\Omega, {\cal F})$. Then\\
 If the entropy $H_Q(P)$ is finite then $P$ is absolutely continuous with respect to $Q$ and 
$$ H_Q(P)=\int log({dP\over dQ})dP.$$
(In particular  if $P$ is not absolutely continuous with respect to $Q$, $H_Q(P)=+\infty$).}\\
This was introduced first by Shannon [ 17 ], for densities: \\
 If $f\in L^1_+(dx)$ is such that $\int_\mathbb R  f(x)dx=1$, Shannon considered $\int_\mathbb R  f(x)logf(x)dx$.\\
 
 Let $\psi$ be the function defined  for $x>0$, by   \begin{eqnarray}\psi(x)=-xlog(x).\end{eqnarray} 
{\bf{Remark 3.2:}}\\
\textit{(1) If $H_Q(P)$ is finite then 
\begin{eqnarray}H_Q(P)=\int log({dP\over dQ})\times {dP\over dQ}dQ=-\int \psi({dP\over dQ})dQ.\end{eqnarray}
In particular,\\\textit{(2) $H_Q(P)=0$ if and only if $P=Q$.}\\
(3) More generally, If $P$ is absolutely continuous with respect to $Q$ then 
\begin{eqnarray}
H_Q(P)<\infty\iff \int log({dP\over dQ})dP<\infty\iff \int (-\psi)({dP\over dQ})dQ<\infty.
\end{eqnarray}}
 \\ The following theorem follows from the monotonicity property (2).\\
{\bf{Theorem B (Dobrushin, [14]):}}  \textit{Let $P,Q$ be probability measures on $(\Omega,\cal F)$, $\cal L$  an algebra of sets 
belonging to $\cal F$, which generates the sigma-algebra $\cal F$, and let $\cal R$ be a family of finite partitions of $\Omega$ whose elements belong to $\cal L$. If every partition consisting of sets from $\cal L$ has a finer partition in $\cal R$, then 
$$H_Q(P)=\sup_{ \Pi\in {\cal R}} S_\Pi (P\mid Q).$$}
* Note that, as remarked by the translator of Pinsker,  in Theorem B of Dobrushin, the condition that the elements of the partitions in $\cal R$ be in $\cal L$ is not necessary.\\

\subsection{A  lemma}

The following lemma will  play an essential role in the rest of the paper.\\
{\bf{Lemma 3.3:}}\\
\textit{Let $(\Omega_i;{\cal F}_i,P_i)$ be a probability space for $i=1,2$, and $\nu$ a probability measure on $(\Omega_1\times \Omega_2,{\cal F}_1\otimes {\cal F}_2)$ with marginals $\nu_1$ on $(\Omega_1,{\cal F}_1)$ and $\nu_2$ on $(\Omega_2,{\cal F}_2)$. Then
\begin{equation}
H_{P_1\times P_2}(\nu)=H_{\nu_1\times\nu_2}(\nu)+H_{P_1}(\nu_1)+H_{P_2}(\nu_2),
\end{equation}
and in particular 
\begin{equation}
H_{P_1\times P_2}(\nu)\ge H_{P_1}(\nu_1)+H_{P_2}(\nu_2),
\end{equation}
 \quad 
\begin{equation}
H_{P_1\times P_2}(\nu)\ge H_{\nu_1\times\nu_2}(\nu).
\end{equation}}\\
{\bf{Proof:}} By Theorem B of Dobrushin,  $H_{P_1\times P_2}(\nu)$ is given by the supremum, over all finite measurable partitions $\Pi_1$ of $\Omega_1$ and $\Pi_2$ of $\Omega_2$, of the sums 
$S_{\Pi_1\times \Pi_2}(\nu\mid P_1\times P_2).$ \\
If $\nu$ is not absolutely continuous with respect to $\nu_1\times \nu_2$, then, by Theorem A, $H_{\nu_1\times \nu_2}(\nu)=+\infty$, and thus 
$$H_{P_1\times P_2}(\nu)\le +\infty=H_{\nu_1\times\nu_2}(\nu)+H_{P_1}(\nu_1)+H_{P_2}(\nu_2).$$
If $\nu$ is absolutely continuous with respect to $\nu_1\times \nu_2$, the equalities  
\begin{eqnarray*}
\nu(E\times F)log{\nu(E\times F)\over P_1\times P_2(E\times F)}= \nu(E\times F)[log{\nu(E\times F)\over \nu_1\times \nu_2(E\times F)}+log {\nu_1(E)\over P_1(E)}+log {\nu_2(F)\over P_2(F)}], E\in \Pi_1,F\in \Pi_2,\\
\end{eqnarray*}
imply 
\begin{eqnarray*}
S_{\Pi_1\times \Pi_2}(\nu\mid P_1\times P_2)=S_{\Pi_1\times \Pi_2}(\nu\mid \nu_1\times \nu_2)+S_{\Pi_1}(\nu_1\mid P_1)+S_{ \Pi_2}(\nu_2\mid  P_2),\hskip 1 cm (E_1)
\end{eqnarray*}
and therefore 
$$S_{\Pi_1\times \Pi_2}(\nu\mid P_1\times P_2)\le H_{\nu_1\times\nu_2}(\nu)+H_{P_1}(\nu_1)+H_{P_2}(\nu_2),$$
from which it  follows  that
$$H_{P_1\times P_2}(\nu)\le H_{\nu_1\times\nu_2}(\nu)+H_{P_1}(\nu_1)+H_{P_2}(\nu_2).$$
   To prove the reverse inequality 
\begin{eqnarray*}H_{P_1\times P_2}(\nu)\ge H_{\nu_1\times\nu_2}(\nu)+H_{P_1}(\nu_1)+H_{P_2}(\nu_2), \hskip 1 cm (E_2)\end{eqnarray*}
consider four arbitrary finite measurable partitions: $ \Pi_1$, $\Delta_1$ of $\Omega_1$, and $\Pi_2$, $\Delta_2$ of $\Omega_2$, and note that we can find a finite measurable partition $\Gamma_i$ of $\Omega_i$ refining both 
$\Pi_i$ and $\Delta_i$, $i=1,2$. Then the partition $\Gamma_1\times \Gamma_2:=\{M\times N:M\in \Gamma_1,N\in \Gamma_2\}$ refines also $\Pi_1\times \Pi_2$. But, in view of inequality $(2)$, we have, for $i=1,2$,
$$S_{\Delta_i}(\nu_i\mid P_i)\le S_{\Gamma_i}(\nu_i\mid P_i).$$
Similarly 
$$S_{\Pi_1\times \Pi_2}(\nu\mid \nu_1\times \nu_2)\le S_{\Gamma_1\times \Gamma_2}(\nu\mid \nu_1\times \nu_2).$$
   So by summing we get, by $(E_1)$ 
   \begin{eqnarray*}S_{\Delta_1}(\nu_1\mid P_1)+S_{\Delta_2}(\nu_2\mid P_2)+S_{\Pi_1\times \Pi_2}(\nu\mid \nu_1\times \nu_2)\le S_{\Gamma_1\times \Gamma_2}(\nu\mid P_1\times P_2))\le H_{P_1\times P_2}(\nu),
   \end{eqnarray*}
 from which $(E_2)$ follows. This proves $(7)$. As trivially  $(7)$ implies $(8)$ and $(9)$, the proof is finished.\\

{\bf{Remark 3.4:}}\\
\textit{The inequality (9) implies that, for fixed $\nu$, the infimum, over all probability measures $P_1$ and $P_2$, of $H_{P_1\times P_2}(\nu)$ is attained for $P_1=\nu_1$ and $P_2=\nu_2$, and the formula $(7)$ shows that it is attained only for these particular values of $P_1$ and $P_2$.}\\

{\bf{Corollary 3.5:}}\\
\textit{Let $\nu$ be a probability measure on a product measurable space. If the entropy $H_{P_1\times P_2}(\nu)$ of $\nu$, with respect to a product probability measure $P_1\times P_2$, is finite , then $\nu$ is absolutely continuous with respect to the product $\nu_1\times \nu_2$ of its marginals, and these marginals are absolutely continuous with respect to $P_1$ and $P_2$ respectively. }\\

\section{Shannon entropy of absolutely continuous processes}
\subsection{Notation}
Let $(\Omega,T,\mu)$ be an invertible ergodic dynamical system. Let ${\cal A}=\{A_1,...,A_k\}$ be a finite partition of $\Omega$. Let $m$ on ${\cal A}^\mathbb Z$ be the product measure with marginals giving equal weights to the atoms of $\cal A$. Let $F=\sum_{j=0}^k a_j1_{A _j}$ be discrete with $a_i\ne a_j$ for $i\ne j$, and $F_n(x)=(F(x),...,F(T^{n-1}x))$, for $x\in \Omega$.  Then if $\mu_n=\mu F_n^{-1}$ and $m_n$ respectively are the restrictions of $\mu$ and $m$ to $\bigvee_{j=0}^{n-1} T^j{\cal A}$, we obtain, with the standard  definition of entropy of a partition 
$$H^\mu(\bigvee_{j=0}^{n-1} T^{-j}{\cal A})=-\int_\Omega  log {d\mu F_n^{-1}\over  dm_n}\circ F_nd\mu.$$

 In this paper  we are interested in the case where $F$ is continuous valued, say real valued, with $(F\circ T^n)$ absolutely continuous (cf. Definition 2.1). 
 In this case,  denoting $l^n$ the Lebesgue measure on $\mathbb R^n$, we  consider the function ${\cal I}_n={\cal I}_n(F)$ defined by:
\begin{eqnarray} {\cal I}_n(x):=-log({d\mu F_n^{-1}\over dl^n})\circ F_n(x), \hskip 0,2 cm x\in \Omega,\end{eqnarray}
together with its integral \begin{equation}
H_n=H_n(F):=\int_\Omega {\cal I}_n d\mu.
\end{equation}
Then, with $\psi$ as in (4), it follows
\begin{eqnarray}
H_n(F)=\int_{\mathbb R^n} \psi({d\mu F_n^{-1}\over dl^n})dl^n.
\end{eqnarray}

We focus on the asymptotic behavior of the sequences ${1\over n}{\cal I}_n$ and ${1\over n}H_n$, which we call, respectively,  the sequences  of \textit{ Shannon information} and \textit{ Shannon entropy } associated to the process $(F\circ T^n)$.\\
Let \begin{eqnarray}
\gamma_0={1\over \sqrt {2\pi}}exp(-{1\over 2}x^2)dl(x),\hskip 2 cm \gamma_n=
\gamma_0^{\otimes n},\hskip 0,12 cm n\ge 1.
\end{eqnarray}

The following quantities are closely related to ${\cal I}_n$ and $H_n$: 
\begin{eqnarray}\hskip 0,5 cm {\cal I}_{n,G}={\cal I}_{n,G}(F)=-log{d\mu F_n^{-1}\over d\gamma_n}\circ F_n \hskip 1 cm (H_{n,G}=H_{n,G}(F)=\int {\cal I}_{n,G}(F)d\mu),\end{eqnarray}
\begin{eqnarray}\hskip 0,5 cm {\cal I}_{n,PM}={\cal I}_{n,PM}(F)=-log{d\mu F_n^{-1}\over d(\mu F^{-1})^{\otimes n}}\circ F_n \hskip 1 cm (H_{n,PM}(F)=\int {\cal I}_{n,PM}(F)d\mu).\end{eqnarray}

The link between these quantities, which behave very much the same, is made precise in the formulas $(46)$ and $(47)$.\\
We shall employ each of the above quantities as seems appropriate. It should be clear that a result formulated using one is easily transformed into a result formulated in terms of the other.\\   
In the case of ergodic $\mathbb Z^2$ action we shall use the following notation\\
 If $F:\Omega \rightarrow \mathbb R$ is  measurable   and $K$ is a finite subset of $\mathbb Z^2$, let $F_K(x)=(F(T^mS^nx))_{(m,n)\in K}$, for $x\in \Omega$.  If the law of $F_K$ is absolutely continuous with respect to $\mid K\mid$-dimensional Lebesgue measure, we denote $f_K$ its density.   In particular if $K=\{(i,j)\in \mathbb Z^2:0\le i\le n-1, 0\le j\le n-1\}$, we denote $F_K$   by $F_{n^2}$ and $f_K$ by $f_{n^2}$. If  the process $(F\circ T^m\circ S^n)_{(m,n)\in \mathbb Z^2}$ is absolutely continuous, as in one dimensional case, we consider \begin{eqnarray}h_n^{(2)}:=-log f_{n^2}\circ F_{n^2},\end{eqnarray} and its integral
\begin{eqnarray}H_n^{(2)}:=\int h_n^{(2)}d\mu.\end{eqnarray}
As our concern is the asymptotic behavior of ${1\over n^2}h_n^{(2)}$, there is no loss of generality if we suppose that $\Omega=\mathbb R^{\mathbb Z^2}$,  $F$ is the projection onto the zero coordinate, $(Tx)_g=x_{g+(0,1)}$, $(Sx)_g=x_{g+(1,0)}$, for $x\in \mathbb R^{\mathbb Z^2}, g\in \mathbb Z^2$, and that $\mu$ is a probability measure whose finite dimensional marginals are absolutely continuous with respect to Lebesgue measure, and which is invariant by the shifts $T$ and $S$. In this case all the densities $f_{n^2}$, for various $n$, will be denoted by $f$ without subscript, so that for every $ x\in \mathbb R^{\mathbb Z^2}$, we have
 \begin{eqnarray}{1\over n^2}h_n^{(2)}(x)=-{1\over n^2}logf(X_{n-1,n-1}^n),\end{eqnarray}

where
$$X^n_{i,j}:=(x_{s,t})_{(s,t)\in I^n_{i,j}},\hskip 3 cm (\alpha_0)$$and $$I^n_{i,j}:=\{(s,j):s=0,...,i\}\cup \{(s,t):0\le s\le n-1,0\le t\le j-1\},\hskip 2 cm (\alpha_1)$$
for $1\le i,j\le n-1$,  $I^n_{0,0}=\{(0,0)\}$ and $I_{i,0}^n=\{(s,0):s=0,...,i\}$.  
Let, for future use,  $$Y^n_{i,j}:=(x_{s,t})_{(s,t)\in I_{i,j}^n,(s,t)\ne (0,0)}.\hskip 3 cm (\alpha_2)$$   In particular
$$X_{n,n}^{n+1}:=(x_{s,t})_{s,t=0,...,n},\hskip 3 cm (a_0)$$and
$$Y_{n,n}^{n+1}:=(x_{s,t})_{s,t=0,...,n, (s,t)\ne 0}.\hskip 3 cm (a_1)$$
Let  $$L:=\{(i,j,n):0\le i,j\le n-1,\hskip 0,2 cm  n\ge 1\}. \hskip 3 cm (a_2)$$ The set inclusion on the  $I^n_{i,j}$, induces a partial order on $L$: we set 
\begin{eqnarray*}
(i,j,n)\le (i',j',n')\iff I^n_{i,j}\subset I^{n'}_{i',j'}\iff (j<j',n\le n') or (j=j',i\le i',n\le n').
\end{eqnarray*}
That is the product of the lexicographical order on $\{(j,i)\} $  and the usual order on $\{n\}$. \\ In the remark below, we single out two properties which we use later:  \\
{\bf{Remark 4.1:}}\\
\textit{ (1) $L$ is directed, and   $C:=\{n-1,n-1,n):n\ge 1\}$ is a cofinal subset of $L$. \\
(2) if $\{(i_k,j_k,n_k):k\in \mathbb N\}$ is an infinite subset of $L$, then there is an infinite  subset $J$ of $\mathbb N$, such that the sequence $((i_k,j_k,n_k))_{k\in J}$ is strictly increasing in $L$, and $\lim_{k\in J} n_k+\infty$. }\\
For $l=(i,j,n)\in L$,  denote ${\cal F}_l$ or $\sigma(X^n_{i,j})$, the sigma-algebra generated by $\{x_{s,t}:(s,t)\in I^n_{i,j}\}$. Then, for $l,l'\in L$, $l\le l'\iff  {\cal F}_l\subset {\cal F}_{l'}$.\\
 If $m$ is a probability measure on $\mathbb R^{\mathbb N^2} $, $m_{i,j}^n$ will denote its restriction to ${\cal F}_l=\sigma(X_{i,j}^n)$,  and then we write 
 $$m_{i,j}^n:=m\mid \sigma(X_{i,j}^n), \hskip 3 cm (r_0)$$
  and in particular $m_{n-1,n-1}^n$ will be denoted simply $m_n$.\\
Two particular probability measures $\pi$ and $\nu$  on $\mathbb R^{\mathbb N^2}$, will be useful for our purpose. Their finite dimensional marginals  $\pi_n=\pi_{n-1,n-1}^n$ and $\nu_n=\nu_{n-1,n-1}^n$ are given by 
\begin{eqnarray}\pi_n=\prod_{(s,t)\in I^n_{n-1,n-1}} f(x_{s,t})d\lambda,\end{eqnarray}
and \begin{eqnarray}\nu_n =f(x_{0,0})\times f(Y_{n-1,n-1}^n)d\lambda,\end{eqnarray}
where $Y_{n-1,n-1}^n$ is as in $(a_1)$ and  $\lambda=\lambda_{n^2}$ denotes $n^2$-dimensional Lebesgue measure.\\
Recall also that \begin{eqnarray}\mu_n=f(X_{n-1,n-1}^n)d\lambda,\end{eqnarray}
where $X_{n-1,n-1}^n$ is as in $(a_0)$.

\subsection
{Convergence of the Shannon  entropy }
We now turn to the dynamical situation. We consider a dynamical system $(\Omega,T,\mu)$ and  denote $AC(\Omega,T,\mu)$ the set of functions $F\in L^2(\mu)$ such that the process $(F\circ T^n)$ is absolutely continuous (i.e. as in Definition 2.1).\\ We establish the convergence of the sequence of Shannon entropy  ${H_n(F)\over n}$  defined by $(11)$, and  give, as in the discrete valued case, an a priori   upper bound  for this limit and a criterion implying that  the process $(F\circ T^n)$ is  Gaussian independent (Corollary 4.5). We also  identify this limit ( Lemma 4.8 ).\\
 In the same way, for $\mathbb Z^2$ action, with $T$ and $S$ as generators, $AC(\Omega,T,S,\mu)$ will denote the set of $F\in L^2(\mu)$ such that the process $(F\circ T^mS^n)_{m,n\in \mathbb Z}$
 is absolutely continuous (i.e. Definition 2.1).\\
The next formula  gives one of the announced links. \\
Recall that $\gamma_n$ denotes the independent Gaussian measure ( cf. (13) ), and $\psi$ is as in (4).\\
{\bf{Lemma 4.2:}}\\
\textit{ Let $F\in AC(\Omega,T,\mu)$. Then
\begin{equation}
{H_n(F)\over n}={1\over n}\int_{\mathbb R^n} \psi({d\mu F_n^{-1}\over d\gamma_n})d\gamma_n+{1\over 2}log(2\pi)+{1\over 2}\mid\mid F\mid\mid_2^2.
\end{equation}}

\textbf{Proof:}
 From formula (12), we get  
\begin{eqnarray*}
H_n
=\int {d\mu F_n^{-1}\over d\gamma_n}\psi({d\gamma_n\over dl^n})dl^n+\int \psi({d\mu F_n^{-1}\over d\gamma_n})d\gamma_n.
\end{eqnarray*}
But, if $h(t):={1\over \sqrt{2\pi}}exp(-{1\over 2}t^2)$ for $t\in \mathbb R$, then 
\begin{eqnarray*}
\int {d\mu F_n^{-1}\over d\gamma_n}\psi({d\gamma_n\over dl^n})dl^n=\int -log{d\gamma_n\over dl^n}d\mu F_n^{-1}
=-\sum_{j=0}^{n-1}  \int log(h(x_j))d\mu F_n^{-1}(x)=-\sum_{j=0}^{n-1} \int log(h(F(T^jx)))d\mu(x)\\
=-n\int log(h(F(x)))d\mu(x)=-n[-{1\over 2}log(2\pi)-{1\over 2}\int F(x)^2d\mu(x)]
=n\times [{1\over 2}log(2\pi)+{1\over 2}\mid\mid F\mid\mid_2^2].\square
\end{eqnarray*}

{\bf{Remark 4.3:}}\\
$(i)$ \textit{The preceding formula $(22)$ can be written as 
\begin{eqnarray}
H_n(F)= -H_{\gamma_n}(\mu F_n^{-1})+{n\over 2}(log 2\pi +\mid\mid F\mid\mid_2^2).\end{eqnarray}}
$(ii)$ \textit{ If $H_n(F) $ is infinite then $H_{n+1}(F)$ is infinite}.\\

In fact, $(ii)$ follows from $(i)$, and $(i)$ from formula $(5)$. \\

{\bf{Lemma 4.4:}}\\
\textit{Let $F\in AC(\Omega,T,\mu)$. Then $(H_n(F))_{n\in N}$ is a sub-additive sequence: for $n,p\in N$
\begin{equation}
H_{n+p}\le H_n+H_p.
\end{equation}
}\\
{\bf{Proof:}}  Inequality (24) follows immediately  from Lemma 3.3 and formula $(23)$.\\

The following corollary is the analogue for the Shannon entropy of the fact that the Kolmogorv entropy of a countable states process is bounded by the entropy of the zerot'h coordinate. The equality case is analogous to the fact that in the Kolmogorov situation, the equality implies that the process is Bernoulli. \\
{\bf{Corollary 4.5:}}\\
\textit{Let $(\Omega,T,\mu)$ be a dynamical system. Then for any $F\in AC(\Omega,T,\mu)$ the sequence $({H_n(F)\over n})$ of Shannon entropies converges to  $Se(F,T)$,  which may be infinite. Moreover:\\
  \begin{eqnarray}Se(F,T)\le {1\over 2}(log(2\pi)+\mid\mid F\mid\mid_2^2),\end{eqnarray} and the equality holds if and only if 
 for every $n$ the law $\mu F_n^{-1}$ of $(F,F\circ T,...,F\circ T^{n-1})$ is the gaussian independent measure $\gamma_n$.} \\

{\bf{Proof:}} The sequence $({H_n\over n}) $ converges to its infimum, since $(H_n)$  is sub-additive by Lemma 4.4.   
 On the other hand,  formula $(23)$  implies  
$ {H_n\over n}\le {1\over 2}(log(2\pi)+\mid\mid F\mid\mid_2^2),$
which yields then the inequality (25). \\
To prove the other statement, note that from formula (23) and  the equality $Se(F,T)=\inf \{{H_n\over n}:n\in N\}$, it follows that the equality in (25) is equivalent to the equalities $H_{\gamma_n}(\mu F_n^{-1})=0,\forall n$. But this is equivalent to the equalities $\mu F_n^{-1}=\gamma_n,\forall n$ by Remark $3.2 \textit{(2)}.\square$\\

  Note that the preceding corollary generalizes to processes a theorem of Shannon that among random variables with fixed variance the maximum of the continuous entropy is achieved by a gaussian variable.\\
A similar proof can yield an $n$ dimensional version of this    theorem of Shannon .\\
 
 {\bf{Proposition 4.6}}\\
 \textit{  Let ${\cal M}_1(n)=\{p\in L_+^1(l^n):\int_{\mathbb R^n} p(x)dl^n(x)=1,\int_{\mathbb R^n} \mid\mid x\mid\mid_2^2p(x)dl^n(x)=n\}$. Then
 $$\sup_{p\in {\cal M}_1(n)} \int_{\mathbb R^n} \psi(p(x))dl^n(x)={n\over 2}(1+log(2\pi)).$$
 Furthermore this supremum is attained  for the gaussian independent density 
 $$p(x)=(2\pi)^{-{n\over 2}}exp(-{1\over 2}(x_0^2+...+x_{n-1}^2))$$
 for $x=(x_0,...,x_{n-1})\in \mathbb R^n$, and this is the only one.}\\

We came to the main definition of this section:\\
{\bf{Definition 4.7:}}\\
\textit{ The Shannon  entropy $Se(F,T)$ of the process $(F\circ T^n)$ is defined by the equality 
\begin{eqnarray}Se(F,T):=\lim_n{H_n(F)\over n}.\end{eqnarray}}
 
 Next we  identify the limit in definition 4.7,  
 the Shannon entropy $Se(F,T)$ of the process $(F\circ T^n)$,  using  conditional entropy, or information ( Lemma 4.8 \textit{(b) (ii)} and \textit{(iii)} ).

{\bf{Lemma 4.8:}}\\
\textit{ (a) Suppose that $H_{n+p}$ is finite. Then \begin{eqnarray}H_{\mu F_n^{-1}\times \mu F_p^{-1}}(\mu F_{n+p}^{-1})=H_n+H_p-H_{n+p},
\end{eqnarray}  
and, for fixed $n$, the sequence $H_{\mu F_n^{-1}\times \mu F_p^{-1}}(\mu F_{n+p}^{-1})$ is increasing in $p$.\\
(b) If $H_n$ is finite for all $n$, then\\
 (i)  $Se(F,T)$ is finite if and only if for every $n$, or for some $n$,  $\sup_p H_{\mu F_n^{-1}\times \mu F_p^{-1}}(\mu F_{n+p}^{-1})$ is finite.\\
 (ii)   \begin{eqnarray}
 \sup_p H_{\mu F_n^{-1}\times \mu F_p^{-1}}(\mu F_{n+p}^{-1})=H_n-n\times Se(F,T).\end{eqnarray}\\
 (iii) If $\nu$ is the law of the process $(F,F\circ T,...)$ and $\nu_n$ is the probability measure on $\mathbb R^\mathbb N$, with  $n+p$-marginal given by $\mu F_n^{-1}\times \mu F_p^{-1}$, for any $p\ge 0$, then \begin{eqnarray*}
H_{\nu_n}(\nu)=H_n-nSe(F,T).
\end{eqnarray*}
In particular, if $Se(F,T)$ is finite, $\nu$ is abolutely continuous with respect to $\nu_n$.}

{\bf{Proof:}} (a) Formula $(27)$ follows from Lemma 3.3 and formula $(23)$. The other property follows from the definition, since, when $\Pi_1$ and $\Pi_2$ are finite partitions of $\mathbb R^n$ and  $\mathbb R^p$, respectively, we have  
$$S_{\Pi_1\times \Pi_2}(\mu F_{n+p}^{-1} \mid \mu F_n^{-1}\times \mu F_p^{-1})\le S_{\Pi_1\times (\Pi_2\times \mathbb R)}(\mu F_{n+p+1}^{-1} \mid \mu F_n^{-1}\times \mu F_{p+1}^{-1}).$$
(b) Put $v_p^n:=H_p-H_{n+p}$, and in particular, for $n=1$, 
 $$u_p:=-v_p^1=H_{p+1}-H_p.$$
 Then, from \textit{(a)} above, the sequence $(u_p)$ is decreasing.
So, as  ${1\over N}\sum_{p=1}^N u_p$ converges to $Se(F,T)$,  $(u_n)$ converges also to $Se(F,T)=\inf_p u_p$. This proves \textit{(i).  (ii)} follows from \textit{(i)} and the equality 
$v_p^n=-u_p-u_{p+1}-...-u_{n+p-1}$. $(iii)$ follows from $(ii).\square$\\

{\bf{Remark 4.9:}}\\
\textit{If $Se(F,T)$ is finite then\begin{eqnarray}
\lim_n {1\over n}\sup_p  H_{\mu F_n^{-1}\times \mu F_p^{-1}}(\mu F_{n+p}^{-1}) =0
\end {eqnarray}
and for fixed $p$,  $H_{\mu F_n^{-1}\times \mu F_p^{-1}}(\mu F_{n+p}^{-1})$  is increasing in $n$, and $z_n:=\sup_p H_{\mu F_n^{-1}\times \mu F_p^{-1}}(\mu F_{n+p}^{-1})$ is  sub-additive and increasing.
}\\
For the convergence of Shannon entropy in the case of $\mathbb Z^2$ action, for $F\in AC(\Omega,T,S,\mu)$, we have  the following \\
{\bf{Lemma 4.10:}}\\
\textit{Let $\mu_n$ and $\pi_n$ be as in $(21)$ and $(19)$ respectively. Then $$\lim_n {1\over n^2}H_n^{(2)}=-\sup_n {1\over n^2}H_{\pi_n}(\mu_n)-\int_\mathbb R f(t)logf(t)dt.$$}\\
{\bf{Proof:}} With notation as in $(a_0)$, 
 we can write $f(X^n_{n-1,n-1})$ in the following form
\begin{eqnarray*}
f(X^n_{n-1,n-1})={f(X^n_{n-1,n-1})\over \prod_{(s,t)\in I^n_{n-1,n-1}} f(x_{s,t})}\times  \prod_{(s,t)\in I^n_{n-1,n-1}} f(x_{s,t}),
\end{eqnarray*}
from which we get

\begin{eqnarray*}
\int log f(X^n_{n-1,n-1})d\mu=\int log {f(X^n_{n-1,n-1})\over \prod_{(s,t)\in I^n_{n-1,n-1}} f(x_{s,t})}d\mu+\sum_{(s,t)\in I^n_{n-1,n-1}} \int log f(x_{s,t})f(x_{s,t})dx_{s,t}.
\end{eqnarray*}
That is,
\begin{eqnarray}
\int log f(X^n_{n-1,n-1})d\mu=H_{\pi_n}(\mu_{n-1,n-1}^n)+n^2\int f(x_{0,0})logf(x_{0,0})dx_{0,0}.\hskip 0,5 cm  (****)
\end{eqnarray}
Put $z_n:=H_{\pi_n}(\mu_{n-1,n-1}^n)$. ( One can see that  $(-z_n)$ is sub-additive.)\\
Let  $R^n_{k,l}=\mathbb Z^2\cap [kn-1,(k+1)n-1]\times [ln-1,(l+1)n-1]$ and denote by $\mu\mid {R^n_{k,l}}$ the restriction of $\mu$ to the coordinates in $R^n_{k,l}$ and similarly for $\pi\mid R^n_{k,l}$. Let $n$ be fixed and $N\ge n$ be an integer. Write $N=p_Nn+r_N=$ where, $p=p_N$, $r=r_N\in \mathbb N$,  $0\le r<n$ and $p\ge 1$. Then clearly we have, by the definition of the conditional entropy
\begin{eqnarray*}
z_N=H_{\pi_N}(\mu^N_{N-1,N-1})\ge H_{\pi_{pn}}(\mu^{pn}_{pn-1,pn-1})=z_{pn}.
\end{eqnarray*}
But,  by the Lemma 3.3  and invariance, we obtain
\begin{eqnarray*}
z_{pn}\ge \sum_{k,l=0}^{p-1} H_{\pi\mid {R^n_{k,l}}}(\mu\mid {R^n_{k,l}})=p^2H_{\pi\mid R^n_{0,0}}(\mu\mid R^n_{0,0})=p^2z_n.
\end{eqnarray*}
So
\begin{eqnarray*}
\limsup_N [ -{1\over N^2}z_N]\le-{z_n\over n^2},
\end{eqnarray*}
which implies 
\begin{eqnarray*}
\limsup_N [ -{1\over N^2}z_N]\le \inf_n [-{z_n\over n^2}].\square
\end{eqnarray*}

As a consequence we can now define an entropy of absolutely continuous process indexed by $\mathbb Z^2$, as follows \\
{\bf{Definition 4.11:}}\\
\textit{The Shannon entropy $Se(F,T,S)$ of the absolutely continuous process $(F\circ T^m\circ S^n)_{(m,n)\in \mathbb Z^2}$ is given by the equality
\begin{eqnarray} Se(F,T,S):=\lim_n [-{1\over n^2}\int logf_{n^2}\circ F_{n^2}d\mu]=\lim_n {1\over n^2}H_n^{(2)}. \hskip 0,4 cm (*)\end{eqnarray}}
Note that here too it follows from  Lemma $4.10$  that
$Se(F,T,S)=-\int f(t)logf(t)dt$ if and only if the process $(F\circ T^m\circ S^n)_{(m,n)\in \mathbb Z^2}$ is independent.\\ 

In order to identify the limit in definition $4.11$, for $F\in AC(\Omega,T,S,\mu)$, we need some further notation.  Let \begin{eqnarray}
g^n_{i,j}:={f(X_{i,j}^n)\over  f(x_{0,0})\times f(Y^n_{i,j})}={d\mu_{i,j}^n\over d\nu_{i,j}^n},
\end{eqnarray}
where $\mu$, $\nu$ are as in (21) and (20), respectively, and  $\mu_{i,j}^n$, $\nu_{i,j}^n$ are as in $(r_0)$.\\
{\bf{Remark 4.12:}}\\
\textit{With $\mu,\nu$, $g_{i,j}^n$ as in $(21),\hskip 0,1 cm (20)$ and $(32)$,  respectively, and $L$ as in $(a_2)$, the following properties  hold:\\
The family $(g^n_{i,j})_{(i,j,n)\in L}$ is a $\nu$ martingale.
\begin{eqnarray}
\int g^n_{i,j}logg^n_{i,j}d\nu\le \int  g^n_{n-1,n-1}log g^n_{n-1,n-1}d\nu.
\end{eqnarray}
 \begin{eqnarray}
\sup_{(i,j,n)\in L} \int g^n_{i,j}logg^n_{i,j}d\nu= \sup_n \int  g^n_{n-1,n-1}log g^n_{n-1,n-1}d\nu.
\end{eqnarray}
 \begin{eqnarray*}
{1\over n^2}h^{(2)}_n=-{1\over n^2}\sum_{i,j=0}^{n-1} log g^n_{i,j} -logf(x_{0,0}).\hskip 2 cm (E)
\end{eqnarray*}}
{\bf{Remark 4.13}}\\
\textit{Notations are as in Remark 4.12. Suppose  $\sup_n \int  g^n_{n-1,n-1}log g^n_{n-1,n-1}d\nu<\infty$. Then if $s_k=(i_k,j_k,n_k)$ is an increasing sequence  in $L$, the martingale  $M_k:=g^{n_k}_{i_k,j_k}$ is uniformly integrable and  converges $\nu$ almost everywhere and in $L^1(\nu)$ to the density $g_s$ of $\mu$ restricted to $ \sigma(\{M_k:k\ge 1\})$ with respect to $\nu$.\\ 
In particular, if $U_n:=g^n_{n-1,n-1}$, then $U_n$ converges $\nu$ almost everywhere to the density $g$ of $\mu$ with respect to $\nu$ on the sigma-algebra $\sigma(\{x_{s,t}:s,t\ge 0\})$ and therefore $g_{i,j}^n=E_\nu(g\mid \sigma(X_{i,j}^n))$. }\\

In fact, by (34),  the hypothesis implies $\sup_{(i,j,n)\in L} \int g^n_{i,j}\mid logg^n_{i,j}\mid d\nu<\infty$,  so that the family $(g^n_{i,j})_{(i,j,n)\in L}$ is uniformly integrable with respect to $\nu$. Here we used the following  \\
{ \bf{Remark 4.14:}}\\
 \textit{
Let $m$ be a finite measure and $\Phi=\{f_i:i\in I\}$ be a family of positive elements in $L^1(m)$ with the property 
\begin{eqnarray*}
\sup_i \int f_ilogf_i dm<\infty.
\end{eqnarray*}
Then (1) $\sup_i \int f_i \mid logf_i\mid dm<\infty.$\\
 (2) The family $\{f_i:i\in I\}$ is uniformly integrable.\\
 (3)  The family $\{(logf_i)^+:i\in I\}$ is uniformly integrable.}\\
 The following lemma  identifies the   Shannon entropy  of $(F\circ S^m\circ T^n)_{(m,n)\in \mathbb Z^2}$, for $\mathbb Z^2$ action (the limit in   Lemma $4.10$).\\
{\bf{Lemma 4.15}}\\
\textit{Let $F\in AC(\Omega,T,S)$, $\mu$ and $\nu$ be as in $(21)$ and $(20)$, respectively. Then \begin{eqnarray}Se(F,T,S)=\lim_n {1\over n^2} H_n^{(2)}=-H_\nu(\mu)-\int f(t)logf(t)dt,\end{eqnarray}
and in particular  $Se(F,T,S)$ is finite if and only if $H_\nu(\mu)<\infty$.}\\
{\bf{Proof:}} 
Equation $(E)$ implies 
\begin{eqnarray*}
{1\over n^2}H^{(2)}_n=-{1\over n^2}\sum_{i,j=0}^{n-1} H_{\nu_{i,j}^n}(\mu_{i,j}^n)-\int f(t)logf(t)dt,
\end{eqnarray*}
so \begin{eqnarray*} {1\over n^2}H^{(2)}_n \ge -H_{\nu^n_{n-1,n-1}}(\mu^n_{n-1,n-1})-\int f(t)logf(t)dt \ge -H_\nu(\mu)-\int f(t)logf(t)dt,  \end{eqnarray*}
 proving one direction. On the other hand, for every $n$, let  $t_n=\int  g^n_{n-1,n-1}log g^n_{n-1,n-1}d\nu=H_{\nu^n_{n-1,n-1}}(\mu^n_{n-1,n-1})$, and $u_n={1\over n^2}\sum_{i,j=0}^{n-1} H_{\nu_{i,j}^n}(\mu_{i,j}^n)$.   If $k> p$,  then
\begin{eqnarray*}
u_k\ge {1\over k^2}\sum_{i=0}^{k-1}\sum_{j=p}^{k-1}H_{\nu_{i,j}^k}(\mu_{i,j}^k)\ge   {1\over k^2}\sum_{i=0}^{k-1}\sum_{j=p}^{k-1}H_{\nu_{p-1,p-1}^p}(\mu_{p-1,p-1}^p)\\
={k(k-p)\over k^2}H_{\nu_{p-1,p-1}^p}(\mu_{p-1,p-1}^p)=(1-{p\over k})H_{\nu_{p-1,p-1}^p}(\mu_{p-1,p-1}^p)=(1-{p\over k})t_p.
\end{eqnarray*}
Hence $\hskip 0,2 cm \liminf_k  u_k\ge t_p,\hskip 0,2 cm$
which implies 
$\hskip 0,2 cm \liminf_k u_k \ge \sup_p t_p, \hskip 0,2 cm$
and  proves the other direction and also the equality
$\hskip 0,2 cm  H_\nu(\mu)=\lim_n u_n, \hskip 0,2 cm$ or that 
$\hskip 0,2 cm {1\over n^2}H^{(2)}_n$ converges to $-H_\nu(\mu)-\int f(t)logf(t)dt.\square$ 

We give examples of non Gaussian and non independent process with finite Shannon entropy.\\
{\bf{Lemma 4.16:}}\\
\textit{Let $(\Omega,\mu,T)$ be a dynamical system. Let $F,G\in L^1(\mu)$. Let $\xi_n:=(F,F\circ T,...,F\circ T^{n-1})$, and $\eta_n:= (G,G\circ T,...,G\circ T^{n-1})$. Suppose that $\xi_n$  is absolutely continuous and that $\xi_n$ and $\eta_n$ are independent.  Then $\xi_n+\eta_n$ is absolutely continuous and 
\begin{eqnarray}
H_n(F+G)\ge H_n(F).
\end{eqnarray}
  }\\
  {\bf{Proof:}} If $H_n(F)=-\infty$, there is nothing to prove. Suppose then that $H_n(F)>-\infty$. Denote by $f$ the density, with respect to $l^n$, of the law $\alpha$ of $\xi_n$, by $\mu$ the law of   $\eta_n$, and by $\nu$ the law of  $\xi_n+\eta_n$. Then $\nu=\alpha\star \mu =gdl^n$, and    $g(t)=\int_{\mathbb R^n} f(t-y)d\mu(y)$, for $l^n$-almost all $t\in \mathbb R^n$. Thus we can write \begin{eqnarray*}
H_n(F+G)=\int_{\mathbb R^n} \psi(g)dl^n=\int_{\mathbb R^n} \psi(\int_{\mathbb R^n}  f(t-y)d\mu(y))dl^n(t)\\
 \ge \int_{\mathbb R^n} [\int_{\mathbb R^n} \psi(f(t-y))d\mu(y)]dl^n(t)=\int_{\mathbb R^n} d\mu(y)\int_{\mathbb R^n} \psi(f(u))dl^n(u)=H_n(F).
\end{eqnarray*}
{\bf{Corollary 4.17:}}\\
\textit{Let $F\in AC(\Omega,T,\mu)$. Assume that $G\in L^1(\mu)$ is such that the two processes $(F\circ T^n)_{n\ge 0}$ and $(G\circ T^n)_{n\ge 0}$ are independent. Then \begin{eqnarray}
Se(F+G,T)\ge Se(F,T).
\end{eqnarray}
If, in addition, $G\in AC(\Omega,T,\mu)$, then\begin{eqnarray*}
Se(F+G,T)\ge \max \{Se(F,T),Se(G,T)\}\ge {1\over 2}(Se(F,T)+Se(G,T)).
\end{eqnarray*}
}
\subsection{Connections between Shannon entropy and Kolmogorov-Sinai entropy}

  We show that if for some $F\in AC(\Omega,T)$, the Shannon entropy $Se(F,T)$ is finite then the Kolmogorv entropy of $(\Omega,T,\mu)$ is infinite (Corollary 4.21). We also  give a new way to describe the Shannon entropy  for $\mathbb Z$ action ( Theorem 4.20 \textit{(b)}), and for $\mathbb Z^2$ action (Remark 4.22) . Finally, we obtain   a criterion of Markovianness (Proposition 4.23). \\ 
The following lemma is important for comparing the Kolmogorov-Sinai entropy to the Shannon entropy.\\
 If $\cal P$ is a finite measurable partition of a probability space $(\Omega,{\cal F},P)$ then the entropy of $\cal P$ will be denoted
  $H({\cal P})$, or $H({\cal P},P)$ or $H^P({\cal P})$.
 
  {\bf{Lemma 4.18:}}\\
\textit{Let $\cal E$ be the set of finite measurable partitions of $\mathbb R$, and $F^{-1}{\cal E}=\{F^{-1}{\cal P}:{\cal P}\in {\cal E}\}$. Then 
$$ H_{\mu F_n^{-1}\times \mu F_p^{-1}}(\mu F_{n+p}^{-1})=\sup_{{\cal B}\in F^{-1}{\cal E}} [H^\mu(\bigvee_{j=0}^{n-1}T^j{\cal B})-H^\mu(\bigvee_{j=0}^{n-1}T^j{\cal B}\mid \bigvee_{j=1}^p T^{-j} {\cal B})].$$}
{\bf{Proof:}} By Theorem B of Dobrushin, 
$$ H_{\mu F_n^{-1}\times \mu F_p^{-1}}(\mu F_{n+p}^{-1})=\sup S_{{\cal P}_1\times {\cal P}_2}(\mu F_{n+p}^{-1}\mid \mu F_n^{-1}\times \mu F_p^{-1}),$$  the supremum being taken over all finite partitions ${\cal P}_1$ of $\mathbb R^n$ and ${\cal P}_2$ of $\mathbb R^p$, of the following particular forms: ${\cal P}_1={\cal P}\times ...\times {\cal P}$, $n$ times, and  ${\cal P}_2={\cal P}\times ...\times {\cal P}$, $p$ times, where ${\cal P}$  run in $\cal E$. But it holds  $$S_{{\cal P}_1\times {\cal P}_2}(\mu F_{n+p}^{-1}\mid \mu F_n^{-1}\times \mu F_p^{-1})=-H({\cal P}_1\times {\cal P}_2,\mu F_{n+p-1}^{-1})+H({\cal P}_1,\mu F_n^{-1})+H({\cal P}_2,\mu F_p^{-1}),$$
  and, for these particular forms of ${\cal P}_i$, $i=1,2$, setting $\hskip 0,1 cm {\cal B}:=F^{-1}{\cal P},\hskip 0,1 cm$ we easily have  \begin{eqnarray*}
H({\cal P}_1,\mu F_n^{-1})=H^\mu (\bigvee_{j=0}^{n-1} T^{-j}{\cal B}),
\end{eqnarray*}
$\hskip 0,1 cm H({\cal P}_2,\mu F_p^{-1})=H^\mu (\bigvee_{j=n}^{n+p-1} T^{-j} {\cal B}),
\hskip 0,1 cm$ 
and 
$\hskip 0,1 cm H({\cal P}_1\times {\cal P}_2,\mu F_{n+p-1}^{-1})=H^\mu (\bigvee_{j=0}^{n+p-1} T^{-j} {\cal B}).$

It follows then that \begin{eqnarray*}
S_{{\cal P}_1\times {\cal P}_2}(\mu F_{n+p}^{-1}\mid \mu F_n^{-1}\times \mu F_p^{-1})
=H^\mu(\bigvee_{j=0}^{n-1}T^j{\cal B})-H^\mu(\bigvee_{j=0}^{n-1}T^j{\cal B}\mid \bigvee_{j=1}^p T^{-j} {\cal B}).\square
\end{eqnarray*}

In the light of lemma 4.18, formula (27) in Lemma 4.8  can be written as
\begin{eqnarray}
\sup_{{\cal B}\in F^{-1}{\cal E}} [ H^\mu(\bigvee_{j=0}^{n-1} T^j{\cal B})-H^\mu(\bigvee_{j=0}^{n-1} T^{j}{\cal B}\mid \bigvee_{j=1}^{p}T^{-j}{\cal B})]=H_n+H_p-H_{n+p}.
\end{eqnarray}
{\bf{Remark 4.19:}}\\
\textit{Formula $(38)$  allows one to obtain conditions which ensure that $H_n$ will be finite.\\ 
One can prove for instance that,  for all $N$,  the following are equivalent:\\
(i) $H_{N+1}$ is finite.\\
(ii) $\sup_{{\cal B}\in F^{-1}{\cal E}}  [H^\mu({\cal B})-H^\mu({\cal B}\mid \bigvee_{j=1}^N T^{-j}{\cal B})]<+\infty$.\\
(ii') $H_{\mu F^{-1}\times \mu F_N^{-1}}(\mu F_{N+1}^{-1})<+\infty$.\\
}\\
The following theorem gives, in particular, the announced new description of the Shannon entropy for $\mathbb Z$ action.\\
{\bf{Theorem 4.20:}}\\
\textit{Let $\cal E$ be as in Lemma 4.18. Suppose $H_n$ finite for all $n$. Then\\
(a)  The following are equivalent:\\
(i) $Se(F,T)$ is finite.\\
(ii) For every  $n$ (or for some $n$),  $\sup_{{\cal B}\in F^{-1} \cal E} [H^\mu(\bigvee_{j=0}^{n-1} T^j{\cal B})-nH^\mu({\cal B}\mid {\cal B}^-)]<+\infty.$\\
(iii)  $\lim_n {1\over n} \sup_{{\cal B}\in F^{-1} \cal E} [H^\mu(\bigvee_{j=0}^{n-1} T^j{\cal B})-nH^\mu({\cal B}\mid {\cal B}^-)]=0$.\\
(b)  The following equality holds 
\begin{eqnarray}
Se(F,T)=H_1+\inf_{{\cal B}\in F^{-1}{\cal E}} [H^\mu({\cal B}\mid {\cal B}^{-})-H^\mu({\cal B})]
\end{eqnarray}}

{\bf{Proof:}}
Formulas $(38)$, $(27)$ and $(28)$ imply the following \begin{eqnarray}
z_n:=\sup_{{\cal B}\in F^{-1}{\cal E}} [H^\mu(\bigvee_{j=0}^{n-1} T^j{\cal B})-H^\mu(\bigvee_{j=0}^{n-1} T^j{\cal B}\mid {\cal B}^-)]=H_n-n\times Se(F,T),
\end{eqnarray}
in which as usual ${\cal B}^-:=T^{-1}{\cal B}\vee T^{-2}{\cal B}\vee...$. . Now 
 \begin{eqnarray*}
H^\mu(\bigvee_{j=0}^{n-1} T^j{\cal B}\mid {\cal B}^-)=nH^\mu({\cal B}\mid {\cal B}^-),\end{eqnarray*}
 so that formula $(40)$ establishes the equivalence between $(i)$ and $(ii)$. As $(i)$ implies $ (iii)$, by Remark 4.9, and trivially $(iii)$ implies $(ii)$, the proof is finished, because taking $n=1$ in formula (40) yields equality (39).\\

 {\bf{Corollary 4.21:}}\\
\textit{Let $(\Omega,T,\mu)$ be an invertible dynamical system. If there exists $F\in AC(\Omega,T,\mu)$ such that the Shannon entropy $Se(F,T)$ of the process $(F\circ T^n)$ is finite then  $(\Omega,T,\mu)$ has infinite entropy.}\\

{\bf{Proof:}} The corollary follows directly from the equivalence between $(i)$ and $(ii)$ in Theorem 4.20.\\

Now, for $\mathbb Z^2$ action,  we establish a formula  analog to  $(39)$ for $\mathbb Z$ action.\\

{\bf{Remark 4.22:}}\\
\textit{Let $\cal E$ be as in Lemma 4.18. Then, for  the stationary absolutely continuous process $(F\circ S^m \circ T^n)$, indexed by $\mathbb Z^2$,  the following holds
 \begin{eqnarray}Se(F,T,S)=-\int f(t)logf(t)dt+\inf_{{\cal P}\in \cal E}  [H^\mu(F^{-1}{\cal P}\mid (F^{-1}{\cal P})^-)- H^\mu(F^{-1}{\cal P})]. \end{eqnarray}
}
In fact, we have $H_\nu(\mu)=\sup_n H_{\nu_n}(\mu_n)$. But
\begin{eqnarray*}
H_{\nu_n}(\mu_n)=\sup S_\Pi (\mu_n\mid \nu_n),
\end{eqnarray*}
the supremum being taken over all partitions $\Pi$ of the form $\Pi=\{\cap_{(i,j)\in I_{n-1,n-1}^n} S^{-i}T^{-j}F^{-1}E: E\in \cal P\}$, where $\cal P$ is a finite partition of $\mathbb R$. But then, with $E_{i,j}=S^{-i}T^{-j}F^{-1}E$, from the definitions of $\mu_n$ and $\nu_n$,  (see (21) and (20) ), it follows 
\begin{eqnarray*}
S_\Pi (\mu_n\mid \nu_n)=\sum_{E\in \cal P} \mu_n(\cap_{(i,j)\in I_{n-1,n-1}^n}E_{i,j})log{\mu_n(\cap_{(i,j)\in I_{n-1,n-1}^n} E_{i,j})\over \nu_n(\cap_{(i,j)\in I_{n-1,n-1}^n} E_{i,j})}=
\\H^\mu(F^{-1}{\cal P})-H^\mu(\bigvee_{(i,j)\in I_{n-1,n-1}^n} S¬^{-i}T^{-j} F^{-1}{\cal P})+H^\mu(\bigvee_{(i,j)\in J_{n-1,n-1}^n} S¬^{-i}T^{-j}F^{-1} {\cal P})=\\
H^\mu(F^{-1}{\cal P})-H^\mu(F^{-1}{\cal P}\mid \bigvee_{(i,j)\in J_{n-1,n-1}^n} S¬^{-i}T^{-j}F^{-1} {\cal P}).
\end{eqnarray*}
Hence 
\begin{eqnarray*}
\sup_n H_{\nu_n}(\mu_n)=\sup_{\cal P}  [H^\mu(F^{-1}{\cal P}) -H^\mu(F^{-1}{\cal P}\mid (F^{-1}{\cal P})^-)],
\end{eqnarray*}
where $(F^{-1}{\cal P})^-=\bigvee_{n\ge 1} \bigvee_{(i,j)\in J_{n-1,n-1}^n} T¬^{-i}S^{-j}F^{-1} {\cal P}).$  

By Lemma $4.15$, we thus obtain (41).\\The following proposition gives a criterion for the Markoviannness of the process  $(F\circ T^n)$.\\
{\bf{Proposition 4.23:}}\\
\textit{Let $F\in AC(\Omega,T,m)$ such that $Se(F,T)$ is finite. Then \\
(i)  The process $(F\circ T^n)$ is Markovian if and only if $Se(F,T)=H_2-H_1.$ More generally\\
(ii) The process $(F\circ T^n)$ has memory $p$ if and only if $Se(F,T)=H_{p+1}-H_p$. }\\

{\bf{Proof:}} We prove $(i)$. The proof of $(ii)$ is similar. If $(F\circ T^n)$ is Markovian, we see by formulas $(13)$ and $(15)$, that $Se(F,T)=H_2-H_1.$\\
For the other direction, put, for $n\ge 1$,   
$$a_n:=\sup_{{\cal A}\in F^{-1}{\cal E}} [H({\cal A})-H({\cal A})\mid T^{-1}{\cal A}\vee...\vee T^{-n}{\cal A}))],$$
 ${\cal R}^n={\cal R}^n({\cal A}):=T^{-2}{\cal A}\vee... \vee T^{-n}{\cal A}$, and  $L=L_n=E[.\mid T^{-1}{\cal A}]-E[.\mid T^{-1}{\cal A}\vee {\cal R}^n]$. Let  $f$ be bounded $F^{-1}{\cal E}$-measurable function, where, as in Lemma 7, ${\cal E}$ denotes the set of finite measurable partitions of $\mathbb R$. We shall show that $L(f)=0$ and the proof will be finished. To do this,  we use a Lemma in [19], according to which for any $\epsilon>0$, there exists $\delta(\epsilon)$, $0<\delta(\epsilon)<\epsilon$, such that for any probability space $(\Omega,{\cal F},\mu)$, and any finite partitions ${\cal P}$ and $\cal Q$ of $\Omega$, the inequality $H^\mu({\cal P})-H^\mu({\cal P}\mid{\cal Q})<\delta(\epsilon)$ implies that $\cal P$ and $\cal Q$ are $\epsilon$-inependent. \\
Recall that following Ornstein, if $\cal P$ and $\cal Q$ are finite measurable  partitions of a probability space $(\Omega,{\cal F},\mu)$, then $\cal P$ is said to be $\epsilon$-independent of  $\cal Q$ if 
$$\sum_{p\in {\cal P}} \mid \mu(p\mid q)-\mu(p)\mid <\epsilon,$$
for all atom $q$ except a set of atoms   of $\cal Q$ which union has a measure less than $\epsilon$.\\  
Fix $\epsilon>0$. Then  there exists ${\cal A}_0$ such that 
$$a_1-\delta(\epsilon)^2 < H({\cal A})-H({\cal A}\mid T^{-1}{\cal A}),$$
for any ${\cal A}$ which is finer than ${\cal A}_0$. 
But, by formula $(15)$, the hypothesis $Se(F,T)=H_2-H_1$ is equivalent to $a_1=a_n$ for any $n$. 
Thus
$$H({\cal A}\mid T^{-1}{\cal A})-H({\cal A}\mid T^{-1}{\cal A}\vee {\cal R}^n({\cal A}))<\delta(\epsilon)^2 .$$
So if we denote, respectively, by $p$, $q$ and $r$, the generic element of ${\cal A}$, $T^{-1}\cal A$ and ${\cal R}^n$, we get
$$\sum_q m(q)[H^{m_q}({\cal A}_q)-H^{m_q}({\cal A}_q\mid {\cal R}^n_q)]<\delta(\epsilon)^2,$$
where $m_q(A)={m(A\cap q)\over m(q)}, \hskip 0,2 cm$ ${\cal A}_q=\{p\cap q:p\in {\cal A}\}$ and similarly for ${\cal R}^n_q$.\\

Let ${\cal Q}_\epsilon:=\{q:H^{m_q}({\cal A}_q)-H^{m_q}({\cal A}_q\mid {\cal R}^n_q)\ge \delta(\epsilon)\}$. 
It follows that 
$$ \sum_{q\in {\cal Q}_\epsilon} m(q) )<\delta(\epsilon), \hskip 1 cm (e_1)$$
and that, for $q\notin {\cal Q}_\epsilon$,  the partitions ${\cal A}_q$ and ${\cal R}^n_q$ are $\epsilon$-independent, under the measure $m_q$, that is, there is $J_q$,  a subfamily of ${\cal R}^n_q$,  such that 
$$\sum_{r\in J_q} m_q(r)>1-\epsilon, \hskip 2 cm (e_2)$$
and
$$\sum_{p\in {\cal A}} \mid m(p\mid q\cap r)-m(p\mid q)\mid<\epsilon, \hskip 1 cm \forall r\in J_q. \hskip 1 cm (e_3)$$
Now we can find ${\cal A}$ finer than ${\cal A}_0$, and $g=g_\epsilon=\sum_{p\in {\cal A}} y_p1_p$ such that $\mid\mid f-g\mid\mid_1<\epsilon$, and $\mid\mid g\mid\mid_\infty\le 2\times \mid\mid f\mid\mid_\infty$. 
Then

$$\mid\mid L(g)\mid\mid_1=\sum_{q,r} \mid \sum_p y_p[m(p\mid q)-m(p\mid q\cap r)]\mid m({q\cap r})=\sum_{q\in {\cal Q}_\epsilon} +\sum_{q\notin {\cal Q}_\epsilon} .$$
In view of $(e_1)$, the first sum in the above equality, is bounded by  $2\mid\mid g\mid\mid_\infty\delta(\epsilon)$.  By $(e_3)$ and $(e_2)$, the second one 
 is bounded by $3\mid\mid g\mid\mid_\infty \times \epsilon$.\\
Therefore 
$$\mid\mid L(g)\mid\mid_1\le 2\mid\mid f\mid\mid_\infty[2\delta(\epsilon)+3\epsilon].$$
  It follows 
$$\mid\mid L(f)\mid\mid_1\le  2\mid\mid f\mid\mid_\infty[2\delta(\epsilon)+3\epsilon]+2\epsilon.$$
This implies $L(f)=0. \square$\\

Note that if $(F\circ T^n)$ is Markovian then  for any $p$, $H_{2p+1}=H_1+2p(H_2-H_1)$
and $H_{2p+2}=H_2+2p(H_2-H_1)$.\\

One might be tempted to introduce an isomorphism invariant 
\textit{  $Se(T)$: 
\begin{eqnarray}
Se(T):=\sup_{F\in AC(\Omega,T\mu)} \{Se(F,T)-H_1(F,T)\}.\end{eqnarray}}\\
It is indeed an invariant, however, it can only take on the 2  values $-\infty$ and $0$.\\

\section{ Entropy rate and Shannon  entropy   }
We shall establish some connections between Shannon  entropy and some concepts developped by Pinsker such as information stability and entropy rate.\\
 We recall some definitions from [14]. First recall that if $Z$ is a random variable then its law is denoted $P_Z$. \\
{\bf{Definition 5.1:}}\\
\textit{Let $\xi=(\xi_n)_{n\ge 1}$ and $\eta=(\eta_n)_{n\ge 1}$ be discrete time stationary  processes.\\
 The entropy rate of $\xi$ with respect to $\eta$ is 
$$\bar{H}_\eta(\xi):=\lim_n{1\over n} H_{P_{(\eta_1,...,\eta_n)}}(P_{(\xi_1,...,\xi_n)}),$$
defined when, for all $j$,  $\xi_j$ and $\eta_j$ take values in the same measurable space,  and when the limit exists.}\\
{\bf{Remark 5.2:}}\\
We can prove, using Lemma 1,   that \textit{ if $\eta$ is the independent Gaussian process, then  for any discrete time real state stationary process $\xi$, the entropy rate of $\xi$ with respect to $\eta$ is well defined.}\\

{\bf{Lemma 5.3:}}\\
\textit{Let $\eta=(X_n)$ be the independent gaussian process, with law $\gamma$ and $\gamma_n$ its projection on the first $n$  coordinates. Let $( \Omega,T,\mu)$ be a dynamical system,  $F\in AC( \Omega,T,\mu)$ and $\nu$ the law of the process $\xi:=(F,F\circ T,...)$. Then\\
(i) $Se(F,T)$ is finite if and only if $\sup_n {1\over n} H_{\gamma_n}(\mu F_n^{-1})<\infty$.\\
(ii) $\nu=\gamma$ if and only if $\lim_n {1\over n}H_{\gamma_n}(\mu F_n^{-1})=0;$\\
(iii) If $H_\gamma(\nu)$ is finite then $\nu=\gamma$.}\\

{\bf{Proof:}} Formula $(23)$  implies 
\begin{eqnarray}
Se(F,T)={1\over 2}(log(2\pi)+\mid\mid F\mid\mid_2^2)-\lim_n {1\over n}H_{\gamma_n}(\mu F_n^{-1}),
\end{eqnarray}
from which we see that $Se(F,T)$ is finite if and only if $\lim_n {1\over n}H_{\gamma_n}(\mu F_n^{-1})$ is finite. 
This proves (i) because by super-additivity \begin{eqnarray*}
\lim_n {1\over n}H_{\gamma_n}(\mu F_n^{-1})=\sup_n {1\over n}H_{\gamma_n}(\mu F_n^{-1}).\end{eqnarray*}\\
Also this last equality together with Remark 3.2 \textit{(2)}  proves (ii).\\
To prove (iii) note that $H_{\gamma_n}(\mu F_n^{-1})\le H_\gamma(\nu)$, and  formula $(23)$ implies thus the inequality
$${1\over n}H_n(F)\ge  {1\over 2}(log(2\pi)+\mid\mid F\mid\mid_2^2)-{1\over n}H_\gamma(\nu)$$
which in turn implies the following one  $$Se(F,T)\ge  {1\over 2}(log(2\pi)+\mid\mid F\mid\mid_2^2).$$
Hence,  by  corollary 4.5,  we have  $\nu=\gamma.\square$\\

{\bf{Corollary 5.4:}}\\
\textit{Let notations be exactly as in Lemma 5.3. Then\\
(i) $Se(F,T)$ is finite if and only if the entropy rate of $\xi$ with respect to $\eta$ is finite.\\ 
(ii) $\nu=\gamma$ if and only if the entropy rate of $\xi$ with respect to $\eta$ vanishes.\\
(iii) \begin{eqnarray}
Se(F,T)={1\over 2}(log(2\pi)+\mid\mid F\mid\mid_2^2)-\bar{H}_\eta(\xi).
\end{eqnarray} }\\

{\bf{Lemma 5.5:}}\\
\textit{ Let $( \Omega,T,\mu)$ be a dynamical system,  $F\in AC( \Omega,T,\mu)$ and $\nu$ the law of the process $\xi:=(F,F\circ T,...)$. Let $P$ be the prduct measure $P:=\mu F^{-1}\otimes \mu F^{-1}\otimes... $ and $P(n)=(\mu F^{-1})^{\otimes n}$ its projection to the first $n$ coordinates. Then\\
(i) If for all $n\in \mathbb N$, $H_n$ is finite (in particular if $Se(F,T)$ is finite)   $\nu$ is locally absolutely continuous with respect to $P$: $\forall n, \mu F_n^{-1}<<(\mu F^{-1})^{\otimes n}$.\\ 
(ii) $Se(F,T)$ is finite if and only if $\sup_n {1\over n}H_{(\mu F^{-1})^{\otimes n}}(\mu F_n^{-1})$ is finite.\\
(iii) $\nu=P\iff \lim_n {1\over n}H_{(\mu F^{-1})^{\otimes n}}(\mu F_n^{-1})=0\iff Se(F,T)=\int_\mathbb R \psi({d\mu F^{-1}\over dl})dl$.\\
(iv) If $H_P(\nu)$ is finite then $\nu=P$.\\
(v) \begin{eqnarray}
Se(F,T)=\int_\mathbb R \psi({d\mu F^{-1}\over dl})dl-\sup_n {1\over n}H_{(\mu F^{-1})^{\otimes n}}(\mu F_n^{-1}).
\end{eqnarray}}.\\

{\bf{Proof:}}  Recalling that ${\cal I}_n(F)$ and  ${\cal I}_{n,PM}(F)$ are defined respectively  by  (10) and (15), we have the following formula
\begin{eqnarray}
{\cal I}_n(F)={\cal I}_{n,PM}(F)-\sum_{j=0}^{n-1} log{d\mu F^{-1}\over dl}\circ F\circ T^j,
\end{eqnarray}
which can be written as
 \begin{eqnarray*}
-{1\over n}log{d\mu F_n^{-1}\over dl^n}\circ F_n(x)=-{1\over n}\sum_{j=0}^{n-1} log{d\mu F^{-1}\over dl}\circ F\circ T^j(x)-{1\over n}log{d\mu F_n^{-1}\over d(\mu F^{-1})^{\otimes n}}(F_n(x)),
\end{eqnarray*}
a proof of which is as follows.\\
First, by taking in formula $(7)$ in Lemma $3.3$, $P_1=(\mu F^{-1})^{\otimes n}$ and $P_2=\mu F^{-1}$, we find \begin{eqnarray*}
H_{(\mu F^{-1})^{\otimes n}\times \mu F^{-1}}(\mu F_{n+1}^{-1})=H_{\mu F_n^{-1}\times \mu F^{-1}}(\mu F_{n+1}^{-1})+H_{(\mu F^{-1})^{\otimes n}}(\mu F_n^{-1})\\+H_{\mu F^{-1}}(\mu F^{-1})
=H_{\mu F_n^{-1}\times \mu F^{-1}}(\mu F_{n+1}^{-1})+H_{(\mu F^{-1})^{\otimes n}}(\mu F_n^{-1}).
\end{eqnarray*}
And next, by taking $p=1$ in formula $(27)$, we find
\begin{eqnarray*}
H_{(\mu F^{-1})^{\otimes n+1}}(\mu F_{n+1}^{-1})=H_n+H_1-H_{n+1}+H_{(\mu F^{-1})^{\otimes n}}(\mu F_n^{-1}),
\end{eqnarray*}
from which it follows,   when $H_{n+1} $ is finite, that  $H_{(\mu F^{-1})^{\otimes n+1}}(\mu F_{n+1}^{-1})$ is finite if and only \\ if $H_{(\mu F^{-1})^{\otimes n}}(\mu F_n^{-1})$ is finite. This, using Theorem A,  implies, by induction on $n$, that if $H_m$ is finite for all $m$, then for every $n$, $\mu F_n^{-1}$ is absolutely continuous with respect to $(\mu F^{-1})^{\otimes n}$. Thus we can write $${d\mu F_n^{-1}\over dl^n}={d\mu F_n^{-1}\over d(\mu F^{-1})^{\otimes n}}\times {d(\mu F^{-1})^{\otimes n}\over dl^n}.$$
So, by taking logarithms, we obtain 
 formula $(46)$, from which follow immediately the formula (45), $(ii)$,  $(iii)$ and $(iv)$.\\

{\bf{Remark 5.6:}} \\
\textit{Formula $(46)$ implies 
\begin{eqnarray*}
H_n=H_{n,PM}+n\times [{1\over 2}log(2\pi)+{1\over 2}\mid\mid F\mid\mid_2^2-H_{\gamma_0}(\mu F^{-1})].
\end{eqnarray*}
In the same way we have\begin{eqnarray}
{\cal I}_n(F)={\cal I}_{n,G}(F)+{n\over 2}log(2\pi)+{1\over 2}\sum_{j=0}^{n-1} F^2\circ T^j
\end{eqnarray}
and therefore (cf. formulas (22) and (23) )
\begin{eqnarray*}
H_n=H_{n,G}+n\times [{1\over 2}log(2\pi)+{1\over 2}\mid\mid F\mid\mid_2^2].
\end{eqnarray*}}

 As a corollary, we obtain the following criterion for independence:\\

{\bf{Corollary 5.7:}}\\
{\it{$Se(F,T)=H_1(F)$ if and only if the process $(F\circ T^n)$ is independent.}}\\
 
  Note that Corollary 5.7   can also be proved  by using formula $(7)$ and $(23)$.  \\
  
  Corollary 5.7  together with formula $(23)$   give the following improvement of Corollary 4.5\\

 {\bf{Corollary 5.8:}}\\
 {\it{ Let $\gamma_0$ be the probability measure with density ${1\over (2\pi)^{1\over 2}}exp(-{1\over 2}x^2)$ with respect to Lebesgue measure $l$ on $\mathbb R$.  Let $F\in AC(\Omega,T,\mu)$. Then\\
 \begin{eqnarray}
Se(F,T)\le -H_{\gamma_0}(\mu F^{-1})+{1\over 2}(log(2\pi)+\mid\mid F\mid\mid_2^2).\end{eqnarray}
 and the equality $$  Se(F,T)=-H_{\gamma_0}(\mu F^{-1})+{1\over 2}(log(2\pi)+\mid\mid F\mid\mid_2^2)$$
 holds if and only if the process $(F\circ T^n)$ is independent.\\
  }}\\
Note, once more, that we see from this corollary that $Se(F,T)= {1\over 2}(log(2\pi)+\mid\mid F\mid\mid_2^2)$ if and only if the process $(F\circ T^n)$ is Gaussian independent [cf. Corollary 4.5].\\
We can also prove the following, which, in particular, improves the inequality (48) in the preceding corollary, and gives a link between the Shannon entropy and information stability.\\

{\bf{Lemma 5.9:}}\\
\textit{If $\xi=(\xi_n)_{n\ge1}$ and $\eta=(\eta_n)_{n\ge 1}$ are discrete time stationary processes, 
the rate of generation of information about $\eta$ by $\xi$ or about $\xi$ by $\eta$ ( following Pinsker) is 
$$\bar{I}(\xi,\eta):=\lim_n {1\over n}H_{P_{(\xi_1,...,\xi_n)}\times P_{(\eta_1,...,\eta_n)}}(P_{(\xi_1,...,\xi_n),(\eta_1,...,\eta_n)}).$$
 The pair $(\xi,\eta)$ is called information stable if $ \bar{I}(\xi,\eta)=0$.\\
Let $(\Omega,T,\mu)$ be an invertible dynamical system and $F\in AC(\Omega,T,\mu)$. Let $\phi$ and $\pi$ be the processes defined by $\phi_n=F\circ T^{-n+1}$, and $\pi_n=F\circ T^n$, for $n=1,2,..$. Then\\
(i)  The pair $(\phi,\pi)$  is information stable if and only if $\lim_n ({H_n\over n}-{H_{2n}\over 2n})=0$. In particular, if $Se(F,T)$ is finite the pair $(\phi,\pi)$ is information stable.\\
(ii)  $Se(F,T)$ is finite if and only if $\sum_{p=0}^\infty {1\over 2^p}
H_{\mu F_{2^p}^{-1}\times \mu F_{2^p}^{-1}}(\mu F_{2^{p+1}}^{-1})<\infty$. Moreover\\
(iii) $Se(F,T)= 
{1\over 2}(log(2\pi)+\mid\mid F\mid\mid_2^2)-H_{\gamma_0}(\mu F^{-1})-{1\over 2} \sum_{p=0}^\infty {1\over 2^p}H_{\mu F_{2^p}^{-1}\times \mu F_{2^p}^{-1}}(\mu F_{2^{p+1}}^{-1}).$}

We also have
 
 {\bf{Remark 5.10:}}\\
 \textit{Le $\Omega,T,\mu)$ be a dynamical system, $F\in AC(\Omega,T,\mu)$ and $\xi$ the process $\xi:=(F,F\circ T,...)$. Then  the entropy rate of $\xi$ with respect to  the independent Gaussian stationary process $\eta$ is given by} 
 \begin{eqnarray*}
 \bar{H}_\eta(\xi)=H_{\gamma_0}(\mu F^{-1})+{1\over 2}\sum_{p=0}^\infty {1\over  2^p}
 H_{\mu F_{2^p}^{-1}\times \mu F_{2^p}^{-1}}(\mu F_{2^{p+1}}^{-1}). \end{eqnarray*}\\

\section{Application to   Gaussian processes }
In this section we express the Shannon entropy $Se(F,T)$ in terms of the spectral measure of   the Gaussian process $X=(X_n)_{n\in \mathbb Z}$,  when $F=X_0$  is the zero
 coordinate function ( Lemma 6.3 ). This enables us (1) to prove that in the class of Gaussian Markovian processes, the Shannon entropy  almost determines the process ( Remark 6.4 ), (2) to show how this entropy changes by linear change of variable ( Corollary 6.5 ), and (3) to prove that all unilateral Gaussian processes with finite Shannon entropy are isomorphic ( Theorem 6.6 ).  \\
 We need first some preliminaries.\\

Let  $(v_n)_{n\ge 0}$ be a stationary sequence of unit vectors in the real Hilbert space $H$, with $(r(n))_{n\in \mathbb Z}$ strictly positive definite sequence (defining $r(-n)=r(n)$), where $r(n)=<v_n,v_0>=<v_{n+k},v_k>$. Let $R_n$ be the $n\times n$ matrix $(R_n)_{ij}=r(i-j)$, for $i,j=0,...n-1$ and $r$ be the vector $r=[r(1),...,r(n-1)]^t$. Then  
 \[
R_n=\left(
\begin{array}{ccc}
  & 1\hskip 0,5 cm r^t  &   \\
  & r\hskip 0,1 cm R_{n-1}  &   \\
  &   &   
\end{array}
\right).
\]

 Evidently for each $n$ there exists a unique vector $a=[a_1,...,a_n]^t$ such that the vector $$w_n:=v_0-\sum_{i=1}^{n-1} a_iv_i$$ is orthogonal to $v_i$ for $i=1,...,n-1$.
Using the orthogonal decomposition $$v_0=w_n+(a_1v_1+...+a_{n-1}v_{n-1}),$$ 
we can prove, by taking scalar products $<v_0,v_0>,...,<v_{n-1},v_0>$, that $a$ is given by the equation
$r=R_{n-1}a$,  or $$a=R_{n-1}^{-1}r.$$

{\bf{Lemma 6.1:}}\\
For any $X^t=(x_0,...,x_{n-1})$ set $Y^t=(x_1,...,x_{n-1})$. Then we have 
\begin{eqnarray}
X^tR_n^{-1}X-Y^tR_{n-1}^{-1}Y={(x_0-\sum_{j=1}^{n-1}a_jx_j)^2\over \mid\mid w_n\mid\mid^2}
\end{eqnarray}
and 
\begin{eqnarray}
det(R_n)= \mid\mid w_n\mid\mid^2 det(R_{n-1}).
\end{eqnarray}

In the same way, note first that we have the equality
\[R_n=
\left(
\begin{array}{ccc}
  & R_{n-1}  \hskip 0,3 cmr^{i}  &   \\
  & ( r^i)^t       \hskip 0,3 cm 1 &   \\
  &   &   
\end{array}
\right)
\]
where $r ^i$ is the vector whose transpose is $(r^i)^t=(r(n-1),...,r(1))$. And for any $n$ there exists a unique vector $b=[b_0,...,b_{n-2}]^t$ such that the vector $u_n$ defined by $$u_n=v_{n-1}-b_0v_0-...-b_{n-2}v_{n-2}$$ is orthogonal to $v_j$ for $j=0,...,n-2$.  \textit{So $v_{n-1}-u_n$ is the projection of $v_{n-1}$ onto the subspace spanned by $v_0,...,v_{n-2}$, and $u_n$ is the projection of $v_{n-1}$ onto the orthogonal of the linear span of $\{v_0,...,v_{n-2}\}$}. Since 
\begin{eqnarray}
v_{n-1}=u_n+(b_0v_0+...+b_{n-2}v_{n-2})
\end{eqnarray}
an "orthogonal decomposition", we obtain, by taking scalar products $<v_0,v_{n-1}>,...,<v_{n-1},v_{n-1}>$,  the equality  $r^i=R_{n-1}b$, or $$b=R_{n-1}^{-1}r^i.$$\\
We have also the following \\

{\bf{Lemma 6.2:}}\\
\textit{If $X=[x_0,...,x_{n-1}]^t$ and $Y=[x_0,...,x_{n-2}]^t$, then 
\begin{eqnarray}
X^tR_n^{-1}X-Y^tR_{n-1}^{-1}Y={(x_{n-1}-\sum_{j=0}^{n-2} b_jx_j)^2\over \mid\mid u_n\mid\mid^2}
\end{eqnarray}
and 
\begin{eqnarray}
det(R_n)=  \mid\mid u_n\mid\mid^2 det(R_{n-1}).
\end{eqnarray}}\\
{\bf{Proof of Lemma 6.2:}} 
 Let \[Q=
\left(
\begin{array}{ccc}
  & I_{n-1} \hskip 0,3 cm 0  &   \\
  &  -b^t\hskip 0,3 cm 1 &   \\
  &   &   
\end{array}
\right)
\]
where $I_{n-1}$ is the identity matrix. Then 
 \[Q^{-1}=
\left(
\begin{array}{ccc}
  & I_{n-1} \hskip 0,3 cm 0  &   \\
  &  b^t\hskip 0,3 cm 1 &   \\
  &   &   
\end{array}
\right).
\]
and we have \[(Q^{-1})^tR_n^{-1}Q^{-1}=
\left(
\begin{array}{ccc}
   & R_{n-1} ^{-1} \hskip 0,3 cm 0 &   \\
  &0...0      \hskip 0,3 cm \alpha &   \\
  &   &   
\end{array}
\right)
\]
where $\alpha={1\over \mid\mid u_n\mid\mid^2}$. In fact this equality is equivalent to 
\[QR_nQ^t=
\left(
\begin{array}{ccc}
   & R_{n-1}  \hskip 0,3 cm 0 &   \\
  &0...0      \hskip 0,3 cm \alpha ^{-1}&   \\
  &   &   
\end{array}
\right)
\]
which can be easily verified.$\square$\\ 

Now let $\Omega=\mathbb R^\mathbb Z$, $\sigma$ the shift transformation, $\mu$ a Gaussian $\sigma$ invariant probability measure determined by a (strictly) positive definite sequence $(r(n))_{n\in \mathbb Z}$, with $r(-n)=r(n)$ for any $n$, so when $r(0)=1$, there exists a probability measure $\nu$ on the unit circle $\mathbb T$ such that $\hat{\nu}(n)=r(n)$.  In this case we shall call $\nu$ the spectral measure. In other words, each $n$ dimensional distribution has a density $\rho_n$ given by \begin{eqnarray}
\rho_n(x_0,...,x_{n-1})={1\over (2\pi)^{n\over 2}(detR_n)^{1\over 2}}\times exp(-{1\over 2}\sum_{i,j=0}^{n-1} (R_n^{-1})_{ij}x_ix_j)
\end{eqnarray}
where $$(R_n)_{ij}=r(i-j)=\int_\Omega x_ix_jd\mu(x), \hskip 0,5 cm i,j=0,...,n-1.$$
 So, if $F(x)=x_0$ for $x\in \Omega$, then, for  $i\ge j$, 
 $$(R_n)_{ij}=\int_\Omega F\circ \sigma^{i-j} Fd\mu.$$ 
 {\bf{Lemma 6.3:}}\\
\textit{Let $P_n$ , $Q_n$ and $Q$ denote the orthogonal projections onto the linear span of $\{X_1,...,X_n\}$ , $ \{X_{-n+1},...,X_{-1}\}$ and  $ \{X_{-1},X_{-2},...\}$ respectively.  Then\\
 a) \begin{eqnarray}Se(F,\sigma)= {1\over 2}log(2\pi)+ 
 log \mid\mid  F-QF\mid\mid_2+{1\over 2}={1\over 2}log(2\pi)+{1\over 2}\int logfd\lambda +{1\over 2}.\hskip 0,2 cm (**)\end{eqnarray}
 where $f$ is the density of $\nu$ with respect to $\lambda$.\\  In particular, if  $\mid\mid F-QF\mid\mid=0$ then $Se(F,\sigma)=-\infty$.\\
 b) If  $\mid\mid F-QF\mid\mid>0$, the following are equivalent\\
(i) The Shannon information  ${1\over n}{\cal I}_n(F)$ converges almost everywhere [resp. in $L^1]$.\\
(ii) ${1\over N}\sum_{j=1}^{N-1} (F-P_jF)^2(\sigma^{N-j})$ converges almost eveywhere  [resp. in $L^1]$.\\
(iii) ${1\over N}\sum_{j=1}^{N-1} (F-P_jF)^2(\sigma^{-j})$ converges almost eveywhere [resp. in $L^1]$.\\
 (iv) ${1\over N}\sum_{j=1}^{N-1} (F-Q_jF)^2(\sigma^j)$ converges almost eveywhere  [resp. in $L^1]$.\\
 (v) ${1\over N}\sum_{n=2}^N log{\rho_n(x)\over \rho_n(\sigma x)}$ converges almost everywhere [resp. in $L^1]$.\\}
{\bf{Proof:}} 
By formula $(54)$ and  lemma 6.1 we get

  \begin{eqnarray}
-log{\rho_n(x_0,...,x_{n-1})\over \rho_{n-1}(x_1,...,x_{n-1})}={1\over 2}log(2\pi)+log\mid\mid F-P_{n-1}F\mid\mid_2+{1\over 2}{(F-P_{n-1}F)^2(x)\over \mid\mid F-P_{n-1}F\mid\mid_2^2}.
\end{eqnarray}

Now, since the process is Gaussian, $P_{n-1}F$ converges to $PF$ almost everywhere and in $L^2$, where $P$ is the projection onto the linear span of $\{F\circ \sigma, F\circ \sigma^2,...\}$, and $PF=E(F\mid \sigma^{-1}{\cal B})$, where $\cal B$ is the Borel sigma-algebra.
So, in the case where $\mid\mid F-PF\mid\mid_2>0$,  it follows from $(56)$, that 
\begin{eqnarray}
\lim_n [-log{\rho_n(x)\over \rho_{n-1}(\sigma x)}]={1\over 2}log(2\pi)+log\mid\mid F-PF\mid\mid_2+{1\over 2}{(F-PF)^2(x)\over \mid\mid F-PF\mid\mid_2^2}.
\end{eqnarray}
Then  the equality\begin{eqnarray}
-{1\over N}\sum_{n=3}^N log{\rho_n(x)\over \rho_{n-1}(\sigma x)}=-{1\over N}\sum_{n=3}^{N-1} log({\rho_n(x)\over \rho_n(\sigma x)})-{1\over N}log\rho_N(x)+{1\over N} log\rho_2(\sigma x)
\end{eqnarray}
 proves that $-{1\over N}log\rho_N(x)$ converges almost surely [ respectively in $L^1$] if and only if $-{1\over N}\sum_{n=3}^{N-1} log({\rho_n(x)\over \rho_n(\sigma x)})$ converges almost surely [respectively in $L^1$].
 Now, by $(56)$ and $(58)$, we obtain the equality $(**)$.\\In the same way, 
we get, by Lemma 6.2
$$-log{\rho_n(x_0,...,x_{n-1})\over \rho_{n-1}(x_0,...,x_{n-2})}={1\over 2}log(2\pi)+{1\over 2}log{det R_n\over det R_{n-1}}+ {1\over 2}{(x_{n-1}-\sum_{j=0}^{n-2} b_jx_j)^2\over \mid\mid u_n\mid\mid^2}.$$

 But, if  $L_{n-2}$ is the orthogonal projection onto the linear span of $X_0,...,X_{n-2}$,  we have
$ \mid\mid u_n\mid\mid^2= \mid\mid  X_{n-1}-L_{n-2} X_{n-1}\mid\mid_2^2= \mid\mid  F-Q_nF\mid\mid_2^2,$
 and thus
\begin{eqnarray*}
{1\over N-1}[-log\rho_N(x_0,...,x_{N-1})+log\rho_1(x_0)]={1\over 2}log(2\pi)+{1\over N-1}\sum_{n=2}^N  
 log \mid\mid  F-Q_nF\mid\mid_2\\
 +{1\over 2(N-1)}\sum_{n=2}^N {(F-Q_nF)^2\circ \sigma^{n-1}(x)\over  \mid\mid  F-Q_nF\mid\mid_2^2}. \hskip 3,5 cm (***)
\end{eqnarray*}\\
Then, in the case where $\mid\mid F-QF\mid\mid>0$, the sequence of Shannon informations ${1\over n}{\cal I}_n(F)$ converges a.e. [ respectively in $L^1$] if and only if ${1\over N}\sum_{j=1}^{N-1} (F-Q_jF)^2(\sigma^j)$ does so.\\
The other statements can be proved in a similar way.$\square$\\

We can see easily from $(**)$ the following\\
{\bf{Remark 6.4:}}\\
\textit{Let $(X_n)$ and $(Y_n)$ be stationary centered Gaussian Markovian processes, with the same $L^2$ norm. Then they have the same Shannon entropy if and only if either they have the same law, or $(Y_n)$ and  $((-1)^nX_n)$ have the same law.}\\

{\bf{Proof:}} Let $ \nu$ be the spectral measure  of $X$ and $\nu'$ be the spectral measure of $Y$. Then if  $P_r(t)=\sum_{n\in \mathbb Z} r^{\mid n\mid}e^{int}$ is the Poisson kernel, we have  $\nu=P_r(t)d\lambda(t)$, for some $r$ and similarly $\nu'=P_{r'}(t)d\lambda(t)$, for some $r'$. On the other hand,  form $(**)$, the equality $Se(X_0,\sigma)=Se(Y_0,\sigma)$ holds if and only if $\mid\mid X_0-QX_0\mid\mid=\mid\mid Y_0-Q'Y_0\mid\mid$, where $Q'$ denotes the projection to the negative coordinates of $Y$. But $QX_0=aX_{-1}$ and similarly $Q'Y_0=bY_{-1}$, for some constants $a,b$. Thus the equality of the respective Shannon entropies is equivalent to $\mid a\mid=\mid b\mid$, or to $a<X_{-1},X_0>=b<Y_{-1},Y_0>$, that is to $ar=br'$. [ Also, one can show by elementary calculus that $\int_\mathbb T logP_r(t)d\lambda(t)=\int_\mathbb T logP_{r'}(t)d\lambda(t)$ if and only if $\mid r\mid =\mid r'\mid$.]$\square$\\

{\bf{Corollary 6.5:}}\\
\textit{Let $(X_n)_{n\in \mathbb Z}$ be Gaussian stationary process with spectral measure $\nu$. Let $g=\sum_{n\in \mathbb Z} a_n e^{int}\in L^2(\nu)$ and $(Y_n)_{n\in \mathbb Z}$ be the stationary Gaussian process such that $Y_0=\sum_{n\in \mathbb Z} a_nX_n$. Then 
$$Se(Y_0,\sigma)=Se(X_0,\sigma)+\int_\mathbb T log(\mid g\mid)d\lambda.$$}
{\bf{Proof:}} By the spectral theorem, if $\nu'$ is the spectral measure of $(Y_n)$, we have $\nu'=\mid g\mid^2\nu$. So, by Szego Theorem $$\mid\mid Y_0-Q'Y_0\mid\mid^2 = exp[\int log(\mid g\mid^2{d\nu\over d\lambda})d\lambda]=exp[\int log(\mid g\mid^2)d\lambda] \mid\mid X_0-QX_0\mid\mid^2.$$
Thus, by $(**)$, we get the result.$\square$\\

Note that, if $g\in L^1(\lambda)$, then $\int_\mathbb T log(\mid g\mid) d\lambda$ is finite if and only if there is $h\in H^1$ such that $\mid g\mid =\mid h\mid$, and in this case, $\int log(\mid g\mid)d\lambda\ge \mid h(0)\mid=\mid \int hd\lambda\mid$.\\
In particular if $a_n=0$ for $n>0$ (or for $n<0$) and $g\in L^1(\lambda)\cap L^2(\nu)$ then  $\int log(\mid g\mid)d\lambda\ge log\mid \int gd\lambda\mid= log\mid  a_0\mid.$\\
More particularly, if  $g\in \mathbb H^1(\mathbb T)$ is an outer function, then $\int log(\mid g\mid)d\lambda= log\mid \int gd\lambda\mid$,  and when, in addition  $\int gd\lambda=0$, we obtain  $Se(Y_0,\sigma)=Se(X_0,\sigma).$\\

It is well known from Ornstein theory that a bilateral gaussian process $X=(X_n)_{n\in \mathbb Z}$ with spectral measure absolutely continuous with respect to Lebesgue measure on the circle is isomorphic to the gaussian independent process. There is interest in considering isomorphism for non-invertible transformations (endomorphisms). The first examples of such isomorphism has been worked out by  Parry [13] and elaborated by Hoffman  and Rudolph [6] ( All the endomorphisms they consider are finite to one.). We consider now the unilateral transformation (endomorphism) associated to gaussian process with spectral measure equivalent to Lebesgue measure $\lambda$. For clarity, if $X=(X_n)_{n\ge 0}$ is a Gaussian process with spectral measure $\nu=fd\lambda$ we consider the endomorphism $T_\nu$ defined on $\mathbb R^\mathbb N$ by $(T_\nu x)_n=x_{n+1}$, for $x\in \mathbb R^\mathbb N$ and $n\ge 0$. Then the shift $T_\nu$ will be isomorphic to the shift $T_\lambda$ if and only if  $log f$ is Lebesgue integrable.\\ To prove  this we recall some useful properties that functions in $\mathbb H^1$ or in $\mathbb H^2$ can have.
First recall that, for $p=1,2$,    $\mathbb H^p$ is the closed subspace of all $f\in L^p(\mathbb T,{d\theta\over 2\pi})$  such that  $\int_{-\pi}^\pi f(t)e^{int}dt=0,n=1,2...$,  and that if $0\le f\in L^1$ then $log f$ is integrable if and only if there is $F\in \mathbb H^2$ such that $f=\mid F\mid^2$ [ [5], Theorem, p.53]. Recall also that an inner function $f$ is an analytic function in the unit disc such that $\mid f(z)\mid \le 1$ and $\mid f(e^{i\theta})\mid =1$ almost everywhere on the unit circle, and an outer function $F$ is an an analytic function in the unit disc of the form \begin{eqnarray*}
F(z)=\alpha exp[{1\over 2\pi}\int_{-\pi}^\pi {e^{i\theta}+z\over e^{i\theta}-z}k(\theta)d\theta]
\end{eqnarray*}

where $k$ is a real-valued integrable function on the circle and $\alpha$ is a complex number with modulus 1 [ [5], p.63,]. For a function $F\in \mathbb H^2$ to be an outer function it is necessary and sufficient that the family $\{z^nF:n=0,1...\}$ span $\mathbb H^2$ [ [5], corollary, p. 101]. Also any non zero function $f\in \mathbb H^1$ can be written in the form $f=gF$ where $g$ is inner and $F$ is outer [ [5], Theorem, p. 63, [4], Theorem 12]. 
\\
In the next theorem the use of Shannon entropy is only to ensure that the logarithm of the density of the spectral measure is integrable.\\

{\bf{Theorem 6.6:}}\\
Let   $\nu$ be a probability measure on the unit circle, equivalent to Lebesgue measure $\lambda$, with density $f$. Then the unilateral shifts $T_\nu$ and $T_\lambda$ are isomorphic if and only if  $Se(X_0,T_\nu)$ is finite, or equivalently $logf$ is Lebesgue integrable. \\

{\bf{Proof:}} Consider the two bilateral gaussian processes $X'$ and $Y'$ with spectral measures $\lambda$ and $fd\lambda$ respectively. Then on the cyclic space $Z_{X_0}=\{X_0\circ T^n:n\in \mathbb Z\}$, $T$ is unitarily equivalent to the multiplication $M_z$ by $z$ on $L^2(\lambda)$. An isomorphism $\phi$ is given by $\phi(X_n)=z^n, n\in \mathbb Z$. Similarily, the same holds for $T$ on the cyclic space $Z_{Y_0}=\{Y_0\circ T^n:n\in \mathbb Z\}$ and the multiplication by $z$ on $L^2(fd\lambda)$, with isomorphism $\psi$: $\psi(Y_n)=z^n, n\in \mathbb Z$. It follows that the action of $T$ on $Z_{X_0}$ is unitarily equivalent to the action of $T$ on $Z_{Y_0}$. Suppose first that $log f$ is integrable. Then  there exists $F\in \mathbb H^2$ such that 
$\sqrt f=\mid F\mid$.  Moreover, there exist an inner function $g$ and an outer function $G\in \mathbb H^2$ such that $F=gG$ and thus $\mid F\mid =\mid G\mid$. Set $x=\phi^{-1} G$, so that $x$ belongs to the closed linear span of $\{X_0,X_1,...\}$, and we have 
\begin{eqnarray*}
<T^nx,x>=<T^n\phi^{-1}G,\phi^{-1}G>=<\phi^{-1}M_z^nG,\phi^{-1}G>\\=<M_z^nG,G>=\int z^n \mid G\mid^2d\lambda=\int z^nfd\lambda=<T^nY_0,Y_0>.
\end{eqnarray*}
On the other hand, if $P=\sum_k a_kz^k$ is a polynomial, the following equalities
\begin{eqnarray*}
\mid\mid X_0-\sum_k a_kT^kx\mid\mid=\mid\mid \phi(X_0)-\sum_k a_k \phi(T^kx)\mid\mid=\mid\mid 1-\sum_k a_kM_z^k\phi(x)\mid\mid\\
=\mid\mid 1-\sum_k a_kM_z^k G\mid\mid=\mid\mid 1-\sum_k a_kz^k G\mid\mid=\mid\mid 1-PG\mid\mid,
\end{eqnarray*}
prove that $X_0$ belongs to the closed linear space generated by $\{T^nx:n\ge 0\}$ if and only if $1$ belongs to the closed linear space generated by $\{z^nG:n\ge 0\}$. But, since $G$ is outer, this later is equal to $\mathbb H^2$ and thus $X_0\in \overline {lin}\{T^nx:n\ge 0\}$. This proves that $T_\nu$ and $T_\lambda$ are isomorphic. The other implication follows from Szeg\"o Theorem.$\square$\\

We end this section with  the following result concerning the speed of convergence in linear prediction:\\

{\bf{Proposition 6.7:}}\\
\textit{Let $\lambda$ be the Lebesgue probability measure on $\mathbb T$, and $\nu$ the spectral measure of a stationary Gaussian process $(X_n)_{n\in \mathbb Z}$. Let $Q$ denote the orthogonal projection onto the closed (in $L^2(\nu)$) linear span of the negative coordinates and $Q_n$ be  the orthogonal projection  onto  the linear span of $\{X_{-n},...,X_{-1}\}$.   Suppose that $X_0\neq QX_0$,  or equivalently $log{d\nu\over d\lambda}$ is Lebesgue integrable. Then:\\ 
  The series 
$\sum_{n=1}^\infty \mid\mid QX_0-Q_nX_0\mid\mid_2^2$ converges if and only if $\nu$ is absolutely continuous and $\nu=e^fd\lambda$, with $\sum_{n=1}^\infty n\mid \hat{f}(n)\mid^2<\infty$}. \\

An equivalent form of  Proposition 3 is \\

{\bf{Remark 6.8:}}\\
\textit{Let $\nu$ be a probability measure on $\mathbb T$. Let $H$, $H_n$ denote the closed subspaces  of $L^2(\nu)$ spanned by $\{e^{ikt}: k\ge 1\}$, and  $\{e^{ikt}: 1\le k\le n\}$ respectively. Let $F$, $F_n$ be the orthogonal projection of the constant function 1 onto $H$ and $H_n$ respectively.  Suppose that 1 is not in $H$. Then the following are equivalent\\
(i) $\sum_{n=1}^\infty  \mid\mid F-F_n\mid\mid^2<\infty$.\\
(ii) $\nu$ is absolutely continuous with respect to the Lebesgue probability measure  $\lambda$ and $\nu=e^fd\lambda$, with $\sum_{n=1}^\infty n\mid \hat{f}(n)\mid^2 <\infty.$}\\

\section{$\mathbb Z^2$ action, pointwise statement}
In this section we consider specifically absolutely continuous  $\mathbb Z^n$ processes, for which we prove pointwise convergence of the Shannon entropy. The case where $n=1$ has already been considered by Barron [1]. However his method can not extend to the higher dimensional case; the idea of our proof is very related to the one by Ornstein and Weiss [10] for the  $\mathbb Z^n$ version of the Shannon Mac Millan Breiman Theorem. The proof is given for $n=2$, but it can be easily generalized.\\
Notations are as in sections 4. Particularly, we refer to (16), (18) for  $h_n^{(2)}$  and to Definition 4.11,  Lemma 4.10 and Lemma 4.15 for $Se(F,T,S)$. Namely, $f_{n^2}$ is the density with respect to Lebesgue measure of the law of $F_{n^2}:=(F\circ T^m\circ S^n)_{m,n=0,...,n-1}$, and $h_n^{(2)}=-logf_{n^2}\circ F_{n^2}$.  The aim of this section is to prove the following theorem:\\
{\bf{Theorem 4.1:}}\\
\textit{Let $T,S$ be commuting measure preserving transformations on the probability space $(\Omega,{\cal F},\mu)$ with ergodic joint action. Let $F\in L^2(\mu)$ such that the process $(F\circ T^m\circ S^n)_{(m,n)\in \mathbb Z^2}$ is absolutely continuous, with law $\nu$. Let  $\nu_0$  be the law of the process  $(F\circ T^m\circ S^n)_{(m,n)\in \mathbb N^2,(m,n)\ne (0,0)}$.\\
Assume that  $Se(F,T,S)$ is finite (which is equivalent to $H_{\mu F^{-1}\times \nu_0}(\nu)<\infty$). Then 
${1\over n^2}h_n^{(2)}$ converges almost everywhere and in $L^1(\mu)$ to $Se(F,T,S)$.\\
In case $Se(F,T,S)=-\infty$ the previous convergence still holds almost everywhere.}\\
 {\bf{Proof:}} We establish first the invariance of $\liminf {1\over n^2}h_n^{(2)}$. Next, with the help of a reduction, we prove that this $\liminf $ is in fact almost everywhere a limit. \\
Recall that, for every $n$,  $X_{n,n}^{n+1}=(x_{s,t})_{s,t=0,...,n}$ is, as in Section 4, formula $(a_0)$.\\
(a): \textit{Let $h_*:=\liminf_n {1\over n^2} h^{(2)}_n$. Then $h_*$ is invariant by each action.}\\
Proof of (a): Let 
$K_{(0,n)\times (0,n)}:=X_{n,n}^{n+1}$,  $K_{(1,n)\times (0,n-1)}:=\{x_{i,j}: 1\le i\le n, 0\le j\le n-1\}$ and  $w_n={1\over n^2}h_n^{(2)}$.   Then  $$w_n\circ S-{ (n+1)^2\over n^2} w_{n+1}=y_n+z_n,$$

where  $\hskip 0,1 cm y_n:={1\over n^2}log {f(K_{(0,n)\times 0,n)})\over \prod_{j=0}^n f(x_{0,j})\times \prod_{i=1}^n f(x_{i,n})\times f(K_{1,n)\times 0,n-1)})}$, and $\hskip 0,1 cm z_n:= {1\over n^2} log( \prod_{j=0}^n f(x_{0,j})\times \prod_{i=1}^n f(x_{i,n}))$.
 Now $\hskip 0,1 cm z_n\hskip 0,1 cm$  converges to 0 almost everywhere by the pointwise ergodic theorems. 
For the first one $\hskip 0,1 cm y_n,\hskip 0,1 cm$  define, for $\hskip 0,1 cm \epsilon>0$,  $\hskip 0,1 cm A_n(\epsilon)=A_n\hskip 0,1 cm $ by 
$$A_n:=\{X: f(X_{n,n}^{n+1})\le e^{-n^2\epsilon}\prod_{j=0}^n f(x_{0,j})\times \prod_{i=1}^n f(x_{i,n})\times  f(K_{(1,n)\times (0,n-1)})\}.$$
Then, a simple calculation yields  $\hskip 0,1 cm \mu(A_n)\le e^{-n^2\epsilon}\hskip 0,1 cm$ and thus  
for $\hskip 0,1 cm \mu\hskip 0,1 cm $ almost all $\hskip 0,1 cm x\hskip 0,1 cm $ there is $\hskip 0,1 cm p\hskip 0,1 cm $ such that  
 $\hskip 0,1 cm -y_n(x)\le \epsilon\hskip 0,1 cm $ for all $\hskip 0,1 cm n\ge p\hskip 0,1 cm $, and this implies 
$\hskip 0,1 cm h_*\le \epsilon +h_*\circ S\hskip 0,1 cm $. Hence 
 $\hskip 0,1 cm h_*\le h_*\circ S$.  It follows that $\hskip 0,1 cm h_*=h_*\circ S$. \\
In the same way, we have also $\hskip 0,1 cm h_*=h_*\circ T\hskip 0,1 cm$ and this  proves (a).\\
We prove now that\\
 (b) \textit{If $ \hskip 0,1 cm \lim_n {1\over n^2}H^{(2)}_n\hskip 0,1 cm$ is finite then the family $\hskip 0,1 cm \{{1\over n^2}h^{(2)}_n:n\ge 1\}\hskip 0,1 cm$ is $\hskip 0,1 cm \mu\hskip 0,1 cm $ uniformly integrable. }\\

Proof of (b):  Recall that $\nu$ and $\mu$ are defined by their respective marginals $\nu_n$,  $\mu_n$ as in (20) and (21) respectively, and  $g_{i,j}^n$ is  as in (32). Also $L$ is as in $(a_2)$ in subsection 4.1. We prove that  the family $\{logg^n_{i,j}:(i,j,n)\in L\}$ is $\mu$ uniformly integrable, and this will imply,  by the equality (E) in Remark 4.12, that $({1\over n^2}h_n^{(2)})$ is $\mu$ uniformly integrable.  For  $l=(i,j,n)$, denote  $g^n_{i,j}$ by $\rho_l$, and let $l_k=(i_k,j_k,n_k)$, be an infinite sequence in $L$. We shall prove that $(log\rho_{l_k})$ contains a weakly convergent subsequence, and this proves (b). By Remark 4.1\textit{(2)},  $(l_k)$ contains a strictly increasing subsequence which we still denote $(l_k)$.  Let ${\cal F}_\infty :=\vee _k {\cal F}_{l_k}$.  By the formula (34) and Remark 4.13, $\rho_{l_k}$ converges $\nu$ almost everywhere to $\rho_\infty:={d\mu_\infty\over d\nu_\infty}$, where $\mu_\infty$ and $\nu_\infty$ are the restrictions of $\mu$ and $\nu$ to ${\cal F}_\infty$ respectively . Also we have
\begin{eqnarray*}
0\ge \int -log\rho_\infty d\mu=-H_{\nu_\infty}(\mu_\infty)\ge -H_\nu(\mu)>-\infty.
\end{eqnarray*}
That is $log\rho_\infty$ is $\mu$ integrable.\\
But  $\sup_k \int \rho_{l_k} log\rho_{l_k}d\mu\le \sup_l \int \rho_l log\rho_l d\mu<\infty$. Hence, by  Remark 4.14, $\{(log\rho_{l_k})^+:k\ge 1\}$ is uniformly integrable with respect to $\mu$. Set $Y_k=log\rho_{l_k}$ and $Y=log\rho_\infty$, so that $Y_k,Y\in L^1(\mu)$, $Y_k$ converges $\mu$ almost everywhere to $Y$ and $\int Y_kd\mu$ converges to $\int Yd\mu$. It follows that $Y_k^+$ converges $\mu$ almost everywhere to $Y^+$, and thus, because $(Y_k^+)$ is $\mu$ uniformly integrable, the convergence holds in $L^1(\mu)$ too. In particular, $\int Y_k^+d\mu$ converges to $\int Y^+d\mu$. So $\int Y_k^-d\mu$ converges to $\int Y^-d\mu$, and thus, since $Y_k^-$ converges $\mu$ almost everywhere to $Y^-$, it converges in $L^1(\mu)$. This proves that $(Y_k)$ converges in $L^1(\mu)$ and a fortiori it is $\mu$ uniformly integrable.\\

Now we proceed to  prove that ${1\over n^2}h_n^{(2)}$ converges $\mu$ almost everywhere. We begin by showing   that\\
(c) \textit{ We can reduce ourselves to the case where the density of the law of the first coordinate is greater than one on its support, and also where $\hskip 0,1 cm \liminf_n {1\over n^2}h_n^{(2)}<0$.}\\
Proof of (c):
Let $F:\mathbb R^{\mathbb Z\times \mathbb Z}\rightarrow \mathbb R$ be the projection to the $(0,0)$ coordinate, with absolutely continuous law with density $f_0={d\mu F^{-1}\over d\lambda}$, and $\alpha>0$. Let  $\phi:\mathbb R\rightarrow \mathbb R$ be the map defined by 
$$\phi(x)={1\over \alpha}\int_{-\infty}^x f_0(t)dt.$$ Put $G=\phi\circ F$. Then  the law of $G$ is absolutely continuous and has density $g_0$ given by $g_0=\alpha1_{[0,{1\over \alpha}]}$.

Let now $\Phi:\mathbb R^{\mathbb Z\times \mathbb Z}\rightarrow \mathbb R^{\mathbb Z\times \mathbb Z}$ be  the map defined by: $(\Phi(x))_{i,j}=\phi(x_{i,j}),\forall i,j\in \mathbb Z$, so that $\hskip 0,1 cm \Phi S^mT^n=S^mT^n\Phi, \hskip 0,1 cm$ and let $\hskip 0,1 cm m:=\mu\circ \Phi^{-1}$. Define $ \theta\hskip 0,1 cm $ by 
 $$\theta(x)=\sup\{t\in \mathbb R: \phi(t)=x\}.$$
Then  the finite dimensional marginals of $m$ are absolutely continuous with respect to Lebesgue measure and  the following relationship holds between the densities $f$ for $\mu$ and $g$ for $m$:
\begin{eqnarray*}
g((u_{i,j})_{i,j=0,...,n-1})=f((\theta(u_{i,j}))_{i,j=0,...,n-1})\prod_{i,j=0,...,n-1}  \theta'(u_{i,j}),
\end{eqnarray*}
with the property: for almost all $ t$, $g_0(t)>0\Rightarrow g_0(t)\ge \alpha.$
It follows  that 
\begin{eqnarray*}
-{1\over n^2}log g((u_{i,j})_{i,j=0,...,n-1})=-{1\over n^2}log f((\theta(u_{i,j}))_{i,j=0,...,n-1})-{1\over n^2}\sum_{i,j=0}^{n-1} log  \theta'(u_{i,j}).
\end{eqnarray*}
But it is easy to see that  $\hskip 0,1 cm \log \theta'\circ F\in L^1(m)\hskip 0,1 cm$ if and only if $\hskip 0,1 cm \int f_0(t)logf_0(t)dt\hskip 0,1 cm$ is finite. Also 
 \begin{eqnarray*}
\int log\theta'\circ Fdm=log\alpha-\int f_0(t)logf_0(t)dt.
\end{eqnarray*} 
So if   $\hskip 0,1 cm \int f_0(t)logf_0(t)dt\hskip 0,1 cm $ is finite, the sequence 
$\hskip 0,1 cm {1\over n^2}\sum_{i,j} log  \theta'(u_{i,j})={1\over n^2}\sum_{i,j} log  \theta'\circ F\circ S^iT^j(u)\hskip 0,1 cm $ converges $\hskip 0,1 cm m\hskip 0,1 cm $ almost everywhere to $\hskip 0,1 cm \int log \theta'\circ Fdm$. Therefore $\hskip 0,1 cm -{1\over n^2}log g((u_{i,j})_{i,j=0,...,n-1})\hskip 0,1 cm $ converges $\hskip 0,1 cm m\hskip 0,1 cm$ almost everywhere if and only if  $\hskip 0,1 cm -{1\over n^2}log f((\theta(u_{i,j}))_{i,j=0,...,n-1})\hskip 0,1 cm$ does so. In this case the corresponding limits (or $\liminf$),  $h_*$ and $g_*$ verifiy
\begin{eqnarray*}
g_*=h_*-log\alpha+\int f_0(t)logf_0(t)dt.
\end{eqnarray*}
So if $h_*$ is $\ge 0$, then if we take $\alpha$ such that 
$$log\alpha> h_*+\int f_0(t)logf_0(t)dt,$$
we obtain $g_*<0$.\\
Therefore, we can suppose that for almost all $t$, if $f_0(t)>0$ then $f_0(t)\ge \alpha\ge 1$ and $h_*:=\liminf {1\over n^2}h_n<0$. \\
This finishes the proof of the announced reduction (c).\\
As in the Ornstein-Weiss case, we prove that\\
(d) \textit{ $\liminf {1\over n^2}h_n^{(2)}$ is almost surely a limit.}\\
Proof of (d):   Put $a:=h_*$. By  (c), we can and do suppose that $a<0$. Let $\epsilon_3>0$, and  $0<\epsilon<\epsilon_3$. Let $l_0\in \mathbb N$, $\epsilon_1,\epsilon_2>0$ and $\delta_1>0$, to be chosen later.\\
First, by the pointwise ergodic theorem, we  find two sequences $(k_l)$ and $(m_l)$ of natural numbers converging to $\infty$ such that $({m_l\over k_l})$ converges to $\infty$ as fast as we wish,  a set $\Omega_\epsilon^1\subset \Omega$, with $\mu(\Omega_\epsilon^1)>1-\delta_1$,  and an integer $N_0(\epsilon_1)$ such that
\begin{eqnarray*}
\forall x\in \Omega_\epsilon^1,\forall N\ge N_0(\epsilon_1), \exists J_N(x)\subset I_{N-1,N-1}^N, \hskip 0,4 cm N^2(1-\delta_1-\epsilon_1)\le card J_N(x),\end{eqnarray*}
and\begin{eqnarray*}
\forall l,\forall (i,j)\in J_N(x),\exists n(i,j,l)\in [k_l,m_l],\hskip 0,4 cm e^{-n^2(i,j,l)(a+\epsilon)}<f((i,j)+X_{n(i,j,l)-1,n(i,j,l)-1}^{n(i,j,l)}).
\end{eqnarray*}
Next, by repeated  uses  of a  Vitali covering type Lemma,  we get an upper estimate of the (Lebesgue) size of the set $\Omega_\epsilon^1$, which enables us to majorize the measure of the set $$\{x:f(X_{N-1,N-1}^N)\le e^{-N^2(a+\epsilon_3)}\}\cap \Omega_\epsilon^1,$$ which will  ensure the convergence of the series $$\sum_N \mu(\{x:f(X_{N-1,N-1}^N)\le e^{-N^2(a+\epsilon_3)}\}\cap \Omega_\epsilon^1),$$ in order to  get $\limsup_n {1\over n^2}h_n^{(2)}\le h_*$. The details are as follows:\\
  
   It is easy to see, from the definition of $a$,  that, given $\epsilon_0>0$, there exist a set $\Omega_\epsilon$, with  $\mu(\Omega_\epsilon )\ge 1-\epsilon_0$, and two strictly increasing sequences of natural numbers $k_l$, $m_l$, with ${m_l\over k_l}$ converging to infinity as fast as one wishes,  such that  
$$\forall l, \forall x\in \Omega_\epsilon, \exists n_l=n_l(x)\in [k_l,m_l] , \hskip 0,3 cm e^{-n_l^2(a+\epsilon)}<f(x)=f(X_{n_l-1,n_l-1}^{n_l}).$$
  Let $J_N(x):=\{(i,j):0\le i,j<N; T^iS^jx\in \Omega_\epsilon\}$.  Then,  by the pointwise ergodic theorem, for any $\delta>0$, there is a measurable set $A_\delta$ with $\mu(A_\delta)>1-\delta$,  and a naural number   $N_0(\epsilon_1,\epsilon,\delta)$,  such that for all $N\ge N_0(\epsilon_1,\epsilon,\delta)$, and for all $x\in A_\delta$, it holds
  $$N^2(\mu(\Omega_\epsilon)-\epsilon_1)< cardJ_N(x)\le N^2(\mu(\Omega_\epsilon)+\epsilon_1).$$ Put $\delta_1=\epsilon_0+\delta$, and $\Omega_\epsilon^1:= \Omega_\epsilon \cap A_\delta$. \\
  It follows $\mu(\Omega_\epsilon^1)\ge 1-\delta_1,$ and  for $x\in\Omega_\epsilon^1$ and $(i,j)\in J_N(x)$ there exists $n(i,j,l)=n(i,j,l)(x) \in [k_l,m_l]$ such that \begin{eqnarray}e^{-n^2(i,j,l)(a+\epsilon)}<f((i,j)+X^{n(i,j,l)}_{n(i,j,l)-1,n(i,j,l)-1}),\end{eqnarray}
 where we  denoted $f(S^iT^jx)$ by  $f((i,j)+X^{n(i,j,l)}_{n(i,j,l)-1,n(i,j,l)-1})$. \\
 
  Note that $R_{i,j,l}(x):=(i,j)+X^{n(i,j,l)}_{n(i,j,l)-1,n(i,j,l)-1}$ is the square with first vertex $(i,j)$ and with side having length $n(i,j,l)$. Then, for any $l$,  $\{R_{i,j,l}(x):(i,j)\in J_N(x)\}$ is a finite cover of $J_N(x)$ by squares. So, due to the freedom, mentioned above, in the choice of $k_l$ and $m_l$, by repeated applications of the  Vitali covering Lemma,[Mattila], for any $\epsilon_2$ there is $N_1(\epsilon_2,\epsilon)$ and  $l_1,...,l_k\ge l_0$,  such that for any $N\ge N_1(\epsilon_2)$, there exist subsets $J_{N,l_1}(x),...,J_{N,l_k}(x)$ of $J_N(x)$, such that  the  squares  $\{R_{i,j,l_s}(x):(i,j)\in J_{N,l_s}(x),s=1...k\}$ are  disjoint, and  there is a subset $J_N^0(x)\subset J_N(x)$ which is covered by  $\{R_{i,j,l_s}(x):(i,j)\in J_{N,l_s}(x),s=1...k\}$ and with 
  $(1-\epsilon_2)card J_N(x)\le card J_N^0(x)$. It follows that, for $N\ge N_0(\epsilon_1,\epsilon,\delta))\vee N_1(\epsilon_2,\epsilon)$ and $x\in \Omega_\epsilon^1$, we have
  $$N^2 u\le card J_N^0(x),$$
  where 
  $\hskip 0,2 cm u:=(1-\epsilon_2)(1-\delta_1-\epsilon_1).\hskip 0,2 cm $
  Then 
  \begin{eqnarray}N^2u\le\sum_{s=1}^k\sum_{(i,j)\in J_{N,l_s}} n^2(i,j,l_s)(x).\end{eqnarray}
 Thus,  since $a<0$, if $\epsilon>0$ is chosen such that  $a+\epsilon< 0$,   
 \begin{eqnarray}
 N^2u(a+\epsilon)\ge\sum_{s=1}^k\sum_{(i,j)\in J_{N,l_s}}  n^2(i,j,l_s)(x) (a+\epsilon).
\end{eqnarray}

 Let $J_N^1(x)$ be the intersection with $I_{N-1,N-1}^N$ of  the union of the squares $\{R_{i,j,l_s}(x):(i,j)\in J_{N,l_s}(x),s=1...k\}$. Then  if $(i,j)\notin J^1_N(x)$ we have that $ 1\le {1\over \alpha}f(x_{i,j})$. 
 In particular $$\alpha\ge 1\Rightarrow 1\le f(x_{i,j}),\forall (i,j)\notin J_N^1(x).$$

  Since $\alpha\ge 1$,  it follows, by (59), that
 \begin{eqnarray*}1\le e^{\sum_{s=1}^k \sum_{(i,j)\in J_{N_,l_s}(x)}  n^2(i,j,l_s)(x) (a+\epsilon)}(\prod_{s=1,...,k,(i,j)\in J_{N,l_s}(x)} f(R_{i,j,l_s}(x)))\times \prod_{(i,j)\notin J_N^1(x)} f(x_{i,j}).
 \end{eqnarray*}
 So by (61), we get for $x\in \Omega_\epsilon^1$,
  \begin{eqnarray}1\le e^{N^2u(a+\epsilon)}(\prod_{s=1,...,k,(i,j)\in J_{N,l_s}(x)} f(R_{i,j,l_s}(x)))\times \prod_{(i,j)\notin J_N^1(x)} f(x_{i,j}).
 \end{eqnarray}
 But the number of all configurations of such disjoint squares is majorised by 
 $C_{N^2}^{[N^2\beta]},\hskip 0,2 cm$ which is majorised by   $ce^{N^2h(\beta,1-\beta)}$,
 where $\beta={1\over q^2}$,  $q$ being the smallest $k_l$'s, and $c$ is a constant. Then by (62), the Lebesgue measure $\lambda(\Omega_\epsilon^1)$ of $\Omega_\epsilon^1$ is, for $N$ big enough, majorized by
 $\hskip 0,1 cm ce^{N^2h(\beta,1-\beta)}\times e^{N^2u(a+\epsilon)}.$
 But easily, 
\begin{eqnarray*}
\mu(\{f(X_{N-1,N-1}^N)\le e^{-N^2(a+\epsilon_3)}\}\cap \Omega^1_\epsilon)\le e^{-N^2(a+\epsilon_3)} \lambda( \Omega^1_\epsilon).
\end{eqnarray*}
So  for $N$ big enough, we then obtain  $\hskip 0,1 cm \mu(\{f(X_{N-1,N-1}^N)\le e^{-N^2(a+\epsilon_3)}\}\cap \Omega^1_\epsilon)\le ce^{-N^2\omega},\hskip 0,1 cm$ where 
$\hskip 0,1 cm \omega=a+\epsilon_3-h(\beta,1-\beta)-u(a+\epsilon).$
But there is a constant $\gamma>0$, such that  for $N$ big enough,  the exponent $\omega$
 is  $>\gamma$. In fact,  if we put $v=\epsilon_1+\delta_1$, then $u=1-v-\epsilon_2(1-v)$, and thus 
 $\hskip 0,1 cm \omega=\epsilon_3-\epsilon-h(\beta,1-\beta)+(a+\epsilon)(v+\epsilon_2(1-v)),\hskip 0,1 cm$
 so that,  for  $0<\gamma<{\epsilon_3-\epsilon\over 2}$, we can choose $\epsilon_1$, $\epsilon_2$, $\delta_1$, $l_0$  and $N_2\ge N_0\vee N_1,$ such that $\forall N\ge N_2$, we have  
the inequality 
$$\hskip 0,1 cm \epsilon_3-\epsilon >2\gamma>\gamma>h(\beta,1-\beta)-(a+\epsilon)(v+\epsilon_2(1-v)).\hskip 0,1 cm$$
Then the series 
$$\sum_N \mu(\{f(X_{N-1,N-1}^N)\le e^{-N^2(a+\epsilon_3)}\}\cap \Omega^1_\epsilon)$$ is convergent. Letting $\delta_1\rightarrow 0$, we get 
$\hskip 0,1 cm \hskip 0,1 cm \limsup_N w_N(x)\le a+\epsilon_3,
\hskip 0,1 cm \hskip 0,1 cm$
and finishes the proof in the case $-\infty <a:=\liminf_N w_N < 0$.\\
If  $\hskip 0,1 cm\liminf w_N=-\infty, \hskip 0,1 cm$ taking $a$ any negative number, the same proof gives $\limsup_N w_N\le a$. \\
This finishes the proof of theorem 4.\\

\centerline{References}
[1] Barron, A.R. : The strong ergodic theorem for densities: Generalized Shannon-McMillan-Breiman theorem. Ann. of Prob., (1985), Vol. 13, No. 4, 1292-1303.\\

[2] Doob, J. L. :  Stochastic processes. Wiley, New York, 1953.\\

[3] Dye, H. A.  : On groups of measure preserving transformations, I. Amer. J. Math., 81, 1959, 119-159.\\

[4] Helson, H. :  Harmonic Analysis. Addison-Wesley (1983).\\

[5] Hoffman, K.:  Banach Spaces of Analytic Functions. Prentice-Hall (1962).\\

[6 ] Hoffman, C., Rudolph, D. J. : Uniform endomorphisms which are isomorphic to Bernoulli shift. Ann. of Math. (2), (2002), 156, No. 1, 79-101.\\

[7] Kolmogorov, A. N.  : Stationary sequences in Hilbert space. Bull. Math. Univ. Moscow 2, no 6  (1941).\\

[8] Mattila P.  : Geometry of Sets and Measures in Euclidean Spaces, Fractals and rectifiability. Cambridge Studies in Advanced Mathematics, 44. Cambridge university Press, Cambridge 1995.\\

[9] Ornstein, D. S. : Ergodic Theory, Randomness, and Dynamical Systems. Yale University Press (1974).\\

[10]  Ornstein, D. S., Weiss B. : The Shannon-Mc Millan-Briman Theorem For A Class Of Amenable Groups. Isr. J. Math. Vol. 44, Mo. 3, 1983, pp. 53-60.\\

 [11]  Parry, W. :  Entropy and Generators in Ergodic Theory. W.A. Benjamin (1969).\\

 [12]  Parry, W.  : Topics in Ergodic Theory. Cambridge University Press (1981).\\

[13 ]  Parry,W.  : Automorphisms of the Bernoulli endomorphism and a class of skew-products. Erg. Th. and Dyn. Syst. 16 (1996), 519-529.\\

 [14] Pinsker,M.S.  :  Information and Information Stability of  Random Variables and Processes. Holden-Day (1964).\\ 

[15] Rosenblatt, M.  : Stationary processes as shifts of functions of independent random variables. J. Math. Mech. 8 (1959) 665-681.\\

[16] Shiriayev, A.N. :  Probability. Second edition, Springer-Verlag (1989).\\

[17]  Shannon, C.E.  : The mathematical theory of communication, Bell System Tech. J. 27 (1948), 379-423; 27 (1948), 623-656.\\

[18] Simon, B. :  The Sharp Form of the Strong Szeg\" o Theorem, To appear in Proc. Conf. on Geometry and Spectral Theory.

[19] Smorodinsky, M.  : Ergodic Theory, Entropy. Lecture Notes, 214, (1971), Springer-Verlag.

[20] Wiener, N. :  Extrapolation, interpolation and smoothing of stationary time series. New York, Wiley (1945).\\

[21] Wiener, N.  : Non-linear problems in random theory. MIT Press, Cambridge, Mass., and Wiley, New York (1958).\\

\end{document}